\documentclass[12pt,letterpaper]{article}
\usepackage{palatino, amssymb,amsmath,amsthm,color,tikz, enumitem}
\usepackage[toc,page]{appendix}
\usepackage{microtype}
\usepackage{url}
\newcommand{\llbracket}{[\![}\newcommand{\rrbracket}{]\!]}

\RequirePackage{graphicx}%

\newcommand{\R}{\mathbb{R}} 

\usepackage[longnamesfirst]{natbib}
\bibpunct[ ]{(}{)}{,}{a}{}{,}

\usepackage[hidelinks]{hyperref}
\providecommand{\doi}[1]{doi: \begingroup\urlstyle{rm}\href{https://doi.org/#1}{\nolinkurl{#1}}\endgroup}

\newtheorem{theorem}{Theorem}[section] %
    \newtheorem{lemma}[theorem]{Lemma}
    \newtheorem{proposition}[theorem]{Proposition}
    \newtheorem{corollary}[theorem]{Corollary}

\theoremstyle{definition} %

\newcommand{\eps}{\varepsilon}
\newcommand{\EE}{\mathbb E}
\newcommand{\PP}{\mathbb P}
\newcommand{\one}{\mathbf 1}

\newcommand{\ed}{\stackrel{d}{=}}

\newcommand{\dR}{{\mathbb R^4_{\uparrow}}}

\definecolor{ao(english)}{rgb}{0.0, 0.5, 0.0}
\makeatletter
\newcommand{\genCal}[1]{\expandafter\newcommand\csname c#1\endcsname{{\mathcal #1}}}
\@for\myi:=A,B,C,D,E,F,G,H,I,J,K,L,M,N,O,P,Q,R,S,T,U,V,W,X,Y,Z\do{%
    \expandafter\genCal\myi
}
\makeatother

\newcommand{\bfa}{\mathbf{a}}
\newcommand{\bfb}{\mathbf{b}}
\newcommand{\bfx}{\mathbf{x}}
\newcommand{\bfy}{\mathbf{y}}
\newcommand{\bfz}{\mathbf{z}}
\newcommand{\bR}{\mathbb{R}}
\newcommand{\bN}{\mathbb{N}}
\newcommand{\bZ}{\mathbb{Z}}

\newcommand{\bc}[1]{\left\{{#1}\right\}}
\newcommand{\br}[1]{\left({#1}\right)}
\newcommand{\bs}[1]{\left[{#1}\right]}
\newcommand{\abs}[1]{\left| {#1} \right|}

\renewcommand{\P}[2][]{\underset{#1}\PP\bs{{#2}}}

\title{The roads taken in the directed landscape}
\author{Pranay Agarwal\and Sourav Sarkar\and B\'alint Vir\'ag}
\date{}

\begin{document}
\maketitle

\begin{abstract}
Every finite-length path in the directed landscape spends almost all of its time on geodesics. More precisely, let \(\mathcal G\) be the union of the interiors of all geodesics. We prove that, almost surely, every path of finite landscape length lies in \(\mathcal G\) outside a set of times of Hausdorff dimension at most \(57/61\). We also construct finite-length paths whose exceptional set has dimension at least \(1/3\). The proof combines quantitative bounds on nearly optimal detours, stochastic domination of the environment outside a geodesic, and a three-geodesic disjointness estimate.
\end{abstract}

\tableofcontents

\section{Introduction}

Every finite-length path in the directed landscape spends almost all of its time on geodesics. Even though the path does not have to maximize length, not even locally, finite length alone forces it onto the geodesic network.

We call a continuous path of finite landscape length a {\bf road}. Let $\mathcal G$ be the union of the interiors of all geodesics. Our main result bounds the exceptional times at which a road leaves $\mathcal G$.

\begin{theorem}\label{t: main}
Almost surely, for every road \(\pi\), the set of times at which \(\pi\) lies outside \(\mathcal{G}\) has Hausdorff dimension at most \(57/61\).
\end{theorem}

\begin{figure}[t]
  \centering
  \includegraphics[width=\linewidth]{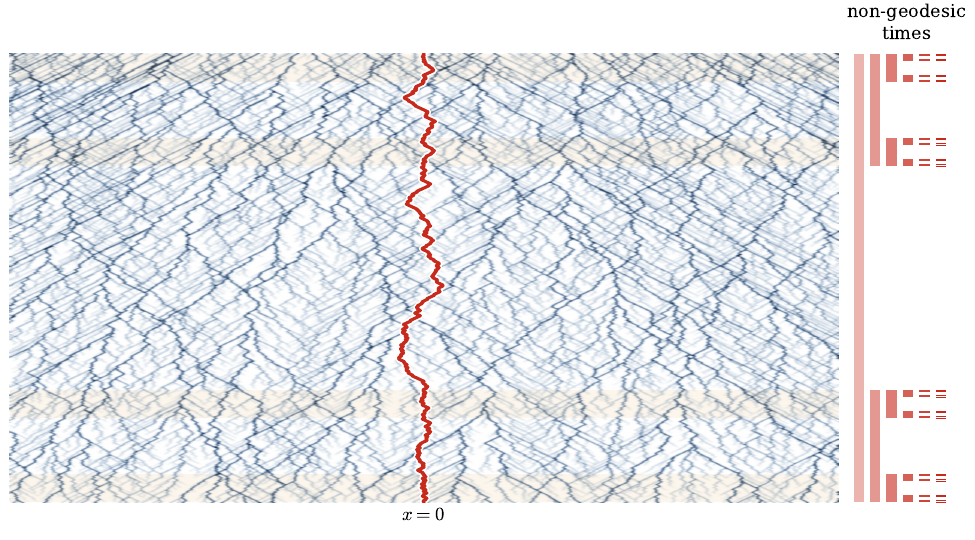}
  \caption{A road winding through a simulated network of geodesics. It runs
  along geodesics almost all the time. It leaves them only on a small fractal
  set of times.}
  \label{fig:intro}
\end{figure}

The exceptional set need not be empty, or even countable. We construct roads with a set of non-geodesic times of dimension $1/3$.

\begin{theorem} \label{t: mainlb}
There exist roads whose set of non-geodesic times has Hausdorff dimension at least \(1/3\) almost surely.
\end{theorem}

Thus the optimal universal bound lies between $1/3$ and $57/61$. 

The directed landscape $\mathcal L$ is a random continuous function on
$\dR=\{(x,s;y,t)\in\mathbb R^4:s<t\}$, constructed by \cite{DOV} as the scaling limit of Brownian last passage percolation. It is conjectured to describe the universal geometry of the KPZ class. Its reverse triangle inequality reads
\[
\mathcal L(p;r)\ge\mathcal L(p;q)+\mathcal L(q;r)
\]
for time-ordered points $p,q,r$. For a continuous path $\pi:[s,t]\to\mathbb R$, its {\bf landscape length} is
\begin{equation}\label{e:length}
|\pi|=\inf\sum_{i=1}^k\mathcal L(\pi(t_{i-1}),t_{i-1};\pi(t_i),t_i),
\end{equation}
where the infimum is over finite subdivisions $s=t_0<\cdots<t_k=t$. We write $\tilde\pi(t)=(\pi(t),t)$. A road has $|\pi|> -\infty$; a geodesic from $u$ to $v$ has length $\mathcal L(u,v)$. The {\bf length deficit} of a path with these endpoints is
\[
\Delta(\pi)=\mathcal L(u,v)-|\pi|.
\]
More precisely, the network in our main theorem is
\begin{equation}\label{e:G-def}
\mathcal G=\big\{(\gamma_{(x,s;y,t)}(r),r):x,y,r,s,t\in\mathbb R,\ s<r<t\big\},
\end{equation}
where $\gamma_{(x,s;y,t)}$ denotes the leftmost geodesic from $(x,s)$ to $(y,t)$.

No fixed deterministic path has finite landscape length almost surely. Our main theorem gives a quantitative form of this distinction: roads cannot spend much time on a prescribed regular path.

\begin{corollary} \label{c: f_int}
For any deterministic locally H\"older $4/61^-$ continuous path \(f\), almost surely, simultaneously for every road \(\pi\), the set of times \(\pi\) intersects \(f\) has Hausdorff dimension at most \(57/61\).
\end{corollary}

Here H\"older $\alpha^-$ means H\"older of every exponent less than $\alpha$. Roads also share the dimension bound of geodesics, despite having weaker pointwise regularity.

\begin{theorem} \label{t: Mink}
Almost surely, every road has Minkowski dimension at most \(4/3\).
\end{theorem}

Geodesics are H\"older $2/3^-$ continuous and have Hausdorff dimension $4/3$; see \cite{dauvergne2022three}. Roads are H\"older $1/2$ continuous, and their modulus can be as bad as this exponent allows. Theorem \ref{t: Mink} follows from a bound on the sum of their squared oscillations. A version of this bound for arbitrary partitions is also an input to the main proof.

\subsection{Detours and the proof of the main theorem}

Leaving a geodesic has a price. The following estimate bounds the probability of a detour that remains disjoint from a geodesic for a prescribed interval while losing very little length.

\begin{theorem}[Detour estimates] \label{t: dmain}
{ Let $0<\delta\le1$. Good-length detours from a geodesic are unlikely.} For a geodesic \(\gamma\) from \((x,r')\) to \((y,t')\) and interval \((r,r+\sigma^3)\subset(r',t')\), the probability that there exists a path \(\pi\) %
with the same endpoints{ , agreeing with \(\gamma\) outside \((r-2\sigma^3,r+3\sigma^3)\)}, disjoint from \(\gamma\) on \((r,r+\sigma^3)\), and length deficit \(|\gamma|-|\pi|\le\delta\sigma\) is at most \(O(\delta\log^{7/3}(1+1/\delta))\). \end{theorem}

We prove this in Section \ref{sec:detour} by separating nearby and far away detours. Nearby detours are controlled by a coupling of the environment outside a geodesic with an independent landscape.

\begin{theorem}[FKG domination outside a geodesic]\label{t:FKG-intro} Let $u\in \dR$.
There exists a coupling of \(\mathcal{L}\) and a copy \(\mathcal{L}'\) such that for every path $\pi$  disjoint from the $u$-geodesic \(\gamma\) we have $|\pi|_{\cL}\le |\pi|_{\cL'}$, and \(\gamma\) is independent of \(\mathcal{L}'\). 
\end{theorem}

Far away detours are controlled by Brownian Gibbs estimates. We apply these at interior mesh times; if a path stays close at every such time, the nearby-detour estimate applies.

The detour bound first compares a road with the geodesic having the same fixed endpoints (Proposition \ref{p:da}). To obtain a comparison uniform over endpoints in spatial intervals, we use the following estimate, proved in Section \ref{sec:triple}.

\begin{theorem}[Triple geodesics] \label{t: 3gd}
For any \(\eta>0\) and \(y>1\), the probability of three disjoint geodesics starting in \([-x,x]\times\{0\}\) and ending in \([-y,y]\times\{1\}\) is at most \(O \br{x^{3/2 -\eta} y^3}\). The exponent $3/2$ is sharp.
\end{theorem}

This gives Proposition \ref{p :cpt_do}: in each time strip, a road can be compared with a geodesic, with a random comparison constant having a polynomial tail. The constants are independent between disjoint time strips, but can be strongly dependent within a strip.

We subdivide time and record which spatial cells a road visits. The partition bounds from Section \ref{sec:regularity} constrain the resulting discrete path. A last-passage estimate then controls how often it can encounter large comparison constants. Combining this count with the sum of the local length deficits gives covers of the non-geodesic times whose $\alpha$-cost tends to zero whenever $\alpha>57/61$. The estimates hold simultaneously for all roads.

\subsection{Related work}

The closest predecessors to our detour estimates are the results of \cite{GangulyHammondNearGround} on near ground states in Brownian last passage percolation. They bound the probability of spatially separated routes with nearly maximal weight and show that slender excursions disjoint from a geodesic incur a substantial weight penalty. These correspond to the far and nearby detour mechanisms used here. Near optimal geodesics (whose weights are close to that of the geodesics) and their local times in exponential directed last passage percolation have been studied in the context of slow-bond TASEP in \cite{SarkarSlyZhang2024}. We formulate our estimates directly in the directed landscape and use them to study all finite-length paths in one environment.

Related excursion arguments enter the study of dynamical last passage percolation. \cite{BhatiaDynamics2026} refines the Ganguly--Hammond excursion argument to establish geodesic-overlap stability and construct a non-trivial continuous stationary reversible dynamics on the directed landscape as a subsequential scaling limit. The companion paper \cite{BhatiaTraces2026} controls geodesic traces through a multiscale analysis of static near-optimal paths. Here the question concerns the geometry of roads in a single landscape.

Disjoint-geodesic estimates and fractal geometry are closely connected. \cite{bates2019hausdorff} study exceptional endpoints supporting disjoint geodesics, while \cite{ganguly2022fractal} construct geodesic local times and determine dimensions of difference profiles. More recently, \cite{GangulySchweinsbergShuai2025} obtain a sharp estimate for three semi-infinite upward geodesics remaining disjoint, with applications to genealogical branch lengths. Their fixed-direction setting differs from the finite-endpoint estimate in Theorem \ref{t: 3gd}.

The geometry of the full geodesic network has been studied through coalescence, branching, and geodesic trees; see \cite{rahman2021infinite,bus22,bhatia2024duality,busani2024stationary}. \cite{dauvergne202527} classify its possible finite geodesic networks. \cite{rahmanVirag2025Trees} reconstruct the landscape from its semi-infinite geodesic trees; their results also yield the countable covering of $\mathcal G$ used here. Other fractal aspects of the landscape are studied in \cite{basu2019fractal,bhatia2022atypical,das2023law}; see also the review \cite{GH24}.

In prelimiting models, geodesic geometry and coalescence have been studied in last passage percolation \cite{Pimentel2013,basu2014last,SeppalainenShen2020,zhang2020optimal,ham20,BSS19, SarkarSlyZhang2024}, polymer models \cite{ras24,alberts2014continuum}, and first passage percolation \cite{Newman1995,dem22}. Results on the environment seen from a geodesic include \cite{MSZ25} for exponential last passage percolation and \cite{dauvergne2022three} for the directed landscape. Geodesics also play a central role in the Brownian map \cite{MQ,AK,LG} and Liouville quantum gravity \cite{gwy21}.

\subsection{Organization of the paper}

Section \ref{sec:prelim} collects preliminary estimates. Section \ref{sec:regularity} proves the regularity and partition bounds for roads. Sections \ref{sec:fkg} and \ref{sec:detour} establish the environment domination and detour estimates, and Section \ref{sec:triple} proves the triple-geodesic bound. The last-passage estimate is proved in Section \ref{sec:last-passage}. Section \ref{sec:trails} combines these ingredients to prove the main theorem and its intersection corollary, and concludes with the construction attaining dimension $1/3$.

\section{Notation and preliminaries}
\label{sec:prelim}

Throughout the article, we  use $c$ to denote positive universal constants that appear during our calculations. Thus, it might change from equation to equation (sometimes in the same display). 

\subsection{Directed landscape estimates}

\begin{theorem}[\cite{DOV}, Corollary~10.7]\label{t:DL-value-bound}
There exists a random constant \(C\) satisfying
\[
\mathbb{P}(C > m) \leq c e^{-d m^{3/2}}
\]
for universal constants \(c, d > 0\) and all \(m > 0\), such that for all \(u = (x, t; y, t + s) \in \mathbb{R}^{4}_\uparrow\), we have
\[
\left| \cL(x, t; y, t + s) + \frac{(x - y)^2}{s} \right| \leq C s^{1/3} \log^{4/3}\left( \frac{2(\|u\| + 2)}{s} \right) \log^{2/3} (\|u\| + 2).
\]
\end{theorem}

The following is a KPZ-scaled version of Theorem 1.2 from \cite{dauvergne2022three}.
\begin{theorem}[Uniform tail bounds for geodesic increments]\label{t:geod-thm-DSV}
	Let $\gamma$ be the directed geodesic from $(0,r_1)$ to $(0,r_2)$. There exists $a>1$, $c > 0$ so that for all $r_1\le s < s + \eps^3=t \le r_2$ we have
	\begin{equation*}
	\EE \bs{a^{|\gamma(s) - \gamma(t)|^3/\eps^6}} < c,  \qquad \EE \bs{a^{|\cL(\gamma(s), s; \gamma(t), t)/\eps|^{3/2}}} < c.
	\end{equation*}
\end{theorem}
This theorem immediately implies a local modulus of continuity bound using Lemma 3.3 \cite{dauvergne2021bulk}, a  Kolmogorov-Chentsov type argument. Together with \cite{dauvergne2021bulk}, Lemma 3.3 we get the following corollary. Note that the local modulus bounds do not get worse as the time span of the whole geodesic gets longer. 
\begin{corollary}
\label{c:local-modulus}
Let $[0,1]\subset [r_1,r_2]$,  let $\gamma$ be the $(0,r_1;0,r_2)$-geodesic. Let $w(t)=\cL(0,r_1;\gamma(t),t)$ { for $t>r_1$, and set $w(r_1)=0$}. Then there is a random constant $C$ satisfying $P(C\ge a)\le e^{1-ca^3}$ so that for all $[t,t+s]\subset [0,1]$ we have
$$
|w(t{ +s})-w({ t})|\le C^2 s^{1/3}|\log(s/2)|^{2/3},\qquad |\gamma(t{ +s})-\gamma({ t})|\le Cs^{2/3}|\log(s/2)|^{1/3}.
$$
\end{corollary}
\subsection{Geodesic network structure of the directed landscape}

\begin{theorem}[Finiteness of $u$-geodesics \cite{dauvergne202527}]\label{t:finite-geodesics}
Almost surely, for every $u\in\dR$, there are only finitely many $u$-geodesics.
\end{theorem}

\begin{theorem}[Absence of bubbles \cite{dauvergne202527,bhatia2024duality}]\label{t:no-bubbles}
Almost surely, for all $(x,s;y,t)\in \dR$ and two $(x,s;y,t)$-geodesics, their set of intersection times is of the form $\{s,t\}\cup I$ for some closed interval $I$. 
\end{theorem}

\subsection{Brownian Gibbs property} The following property enjoyed by the parabolic Airy line ensemble $\cA$ was established in \cite[Theorem $3.1$]{CH}. For any $k,\ell\in \{0,1,2,\ldots\}$ and $a<b\in \R$, let $\cF$ denote the sigma algebra generated by the set
 \[\{\cA_i(x): (i,x)\notin\{k+1,k+2,\ldots,k+\ell\}\times(a,b)\}\,.\]
 Then the conditional distribution of $\cA_{\big|\{k+1,k+2,\ldots,k+\ell\}\times[a,b]}$ given $\cF$ is equal to the law of $\ell$ independent Brownian bridges $B_1,B_2,\ldots, B_\ell:[a,b]\mapsto \R$ with diffusion parameter $2$ such that $B_i(a)=\cA_{k+i}(a)$ and $B_i(b)=\cA_{k+i}(b)$, for all $i=1,2,\ldots, \ell$, conditioned on the event
 \[\cA_k(r)>B_1(r)>B_2(r)>\ldots>B_\ell(r)>\cA_{k+\ell+1}(r)\]
for all $r\in [a,b]$. For this to hold for $k=0$, we set $\cA_0\equiv\infty$.

There are other formulations of the Brownian Gibbs property for the Airy line ensemble. The following formulation from \cite{dauvergne2024wiener} shall be especially useful for proving the triple geodesic exponent (Section \ref{sec:triple}).
\begin{theorem}{\cite[Theorem 3.8]{dauvergne2024wiener}} \label{t: weiner}
Fix $t \ge 1$, $k,\ell \in \mathbb{N}$ and $a \in \mathbb{R}^\ell$ such that $
a_j + t < a_{j+1} \quad \text{for all } j \in \llbracket 1,\ell-1 \rrbracket.$
Define $W(a) = \bigcup_{j=1}^{\ell} \big( \llbracket 1,k \rrbracket \times (a_j, a_j+t) \big).$

Then there exists a random sequence of continuous functions $
\mathcal{L}^a = \mathcal{L}^{t,k,a} = \{ \mathcal{L}_i^a : \mathbb{R} \to \mathbb{R},\ i \in \mathbb{N} \}$ such that the following points hold:

\begin{enumerate}
\item Almost surely, $\mathcal{L}^a$ satisfies
\[
\mathcal{L}_i^a(r) > \mathcal{L}_{i+1}^a(r)
\quad \text{for all pairs } (i,r) \notin W(a).
\]

\item The ensemble $\mathcal{L}^a$ has the following Gibbs property.  
Let $m \in \mathbb{N}$ and $a < b \in \mathbb{R}$ and set $S = \llbracket 1,m \rrbracket \times [a,b].$ Then conditional on the values of $\mathcal{L}_i^a(r)$ for $(i,r) \notin S$, the distribution of
$\mathcal{L}_i^a(r)$, $(i,r) \in S$, is given by $m$ independent Brownian bridges $B_1,\dots,B_m$ from $(a,\mathcal{L}_i^a(a))$ to $(b,\mathcal{L}_i^a(b))$ for $i \in \llbracket 1,m \rrbracket$, conditioned on the event
\[
B_i(r) > B_{i+1}(r) \quad \text{whenever } (i,r) \notin W(a).
\]
\item Let $\mathrm{ORD}$ denote sequences of continuous functions $(f_i)$ such that $f_i > f_{i+1}$ for all $i \in \bN$, then we have that
\begin{equation}
\mathbb{P}\big( \mathcal{L}^a  \in \mathrm{ORD} \big)
\ge
\exp\!\big( - e^{O_k(1+\ell \log \ell)}\, t^3 \big),
\end{equation}
and conditional on the event $\{ \mathcal{L}^a \in \mathrm{ORD} \}$, the ensemble $\mathcal{L}^a$ is equal in law to the parabolic Airy line ensemble $\cA$.
\end{enumerate}
\end{theorem}

\subsection{Extended Airy Sheet} \label{ss: airy}

To study multipath geometry within the directed landscape $\mathcal{L}$, we utilize the framework of the extended Airy sheet constructed in \cite{dauvergne2022disjoint}. We introduce some notation to make the upcoming discussion easier to follow. Let $\bR^k_> = \bc{\bfx \in \bR^k \;|\; x_1 >\cdots > x_k},$ similarly for $\leq, \geq,$ and $<$.

Recall that the classical Airy sheet $\mathcal{S}\colon \mathbb{R}^2 \to \mathbb{R}$ describes the last-passage value (or landscape increment) over a unit time horizon. The \emph{extended Airy sheet}, denoted by $\mathcal{S}(\mathbf{x}, \mathbf{y})$ for vector-valued inputs $\mathbf{x} = (x_1, \dots, x_k)$ and $\mathbf{y} = (y_1, \dots, y_k)$ in $\mathbb{R}^k_\le$, generalizes this construction by evaluating the maximal weight of $k$ mutually disjoint directed paths (or $k$-optimizers) connecting the source vector $\mathbf{x}$ at time $0$ to the target vector $\mathbf{y}$ at time $1$. Specifically, $\mathcal{S}(\mathbf{x}, \mathbf{y})$ corresponds to the restriction of the extended directed landscape $\mathcal{L}$ to multiple disjoint paths:
\begin{equation}
    \mathcal{S}(\mathbf{x}, \mathbf{y}) = \mathcal{L}(\mathbf{x}, 0; \mathbf{y}, 1).
\end{equation}
Continuity of the extended Airy sheet will be of particular interest. To understand this property it will be useful to look at the parabolically shifted sheet $\mathcal{R}$, defined as:
\begin{equation}
    \mathcal{R}(\mathbf{x}, \mathbf{y}) := \mathcal{S}(\mathbf{x}, \mathbf{y}) + \sum_{i=1}^{k} (x_i - y_i)^2.
\end{equation}
In particular, we have that
\begin{lemma} [{\cite[Lemma 5.6]{dauvergne2022disjoint}}] \label{l: esc}
Take any $k \in \mathbb{N}$, $\mathbf{u} = (\mathbf{x}, \mathbf{y})$, $\mathbf{u}' = (\mathbf{x}', \mathbf{y}') \in \mathbb{R}_{\le}^k$ with $\|\mathbf{u} - \mathbf{u}'\|_2 < 1$, and $a > 0$. Then
\begin{equation*}
    \mathbb{P}\left(\left|\mathcal{R}(\mathbf{x}', \mathbf{y}') - \mathcal{R}(\mathbf{x}, \mathbf{y})\right| > a \sqrt{\|\mathbf{u} - \mathbf{u}'\|_2}\right) < c e^{-d a^{3/2}},
\end{equation*}
for some constants $c, d > 0$ depending only on $k$ .
\end{lemma}
We again use the above lemma along with \cite[Lemma 3.3]{dauvergne2021bulk} to get a local modulus of continuity bound.
\begin{corollary}\label{c: esc-unif}
Fix $k\in\bN$. There exist constants $c,d>0$ depending only on $k$ such that for every $\bfz\in\bR^{2k}$ there is a random constant $C_{\bfz}$ with $\P{C_{\bfz}>a}\le c\,e^{-da^{3/2}}$ for all $a>0$, such that for all $\mathbf u=(\bfx,\bfy)$, $\mathbf u'=(\bfx',\bfy')\in \bfz+[-1/2,1/2]^{2k}$,
\[
\abs{\cR(\bfx',\bfy')-\cR(\bfx,\bfy)}\le C_{\bfz}\,\|\mathbf u-\mathbf u'\|_2^{1/2}\log^{2/3}\!\Big(\frac{c}{\|\mathbf u-\mathbf u'\|_2}\Big).
\]
\end{corollary}

\subsection{Coupling results}

\begin{theorem}[Strassen's Theorem, \cite{strassen1965}, \cite{lindvall1999}]\label{t:strassen}
Let \(X\) and \(Y\) be random elements taking values in a partially ordered Polish space \((S,\preceq)\) with closed order relation. The following are equivalent:
\begin{enumerate}
\item \(\mathbb{E}[f(X)]\le \mathbb{E}[f(Y)]\) for every bounded measurable non-decreasing function \(f:S\to\R\);
\item there exists a coupling \((X',Y')\) of \(X\) and \(Y\) such that \(X'\preceq Y'\) almost surely.
\end{enumerate}
\end{theorem}

\begin{lemma}[Gluing Lemma, see \cite{villani2008optimal}]\label{l:gluing}
Let \(X_1,X_2,X_3\) be random variables taking values in Polish spaces. Suppose we are given a coupling of the pair \((X_1,X_2)\) and a coupling of the pair \((X_2,X_3)\) (possibly defined on different probability spaces). Then there exists a joint coupling of the triple \((X_1,X_2,X_3)\) that realizes both given couplings as its marginals.
\end{lemma}

\subsection{Hausdorff dimension}

For a set $A \subseteq \mathbb{R}$ and $s \geq 0$, the $s$-dimensional Hausdorff measure of $A$ is defined as
\[
    \cH^{s}(A) := \lim_{\delta \downarrow 0} \cH^{s}_\delta(A), \qquad \cH^{s}_\delta(A) = \inf \bc{ \sum_{i=1}^\infty \abs{U_i}^s : A \subset \bigcup_{i=1}^\infty U_i, \abs{U_i} < \delta}.
\]
where the infimum is taken over all countable covers $\{U_i\}$ of $A$ by Borel sets, and $|U_i|$ denotes the diameter of $U_i$.

The Hausdorff dimension of $A$ is defined as
\[
\dim_H(A) = \inf \{ s \geq 0 : \cH^s(A) = 0 \}.
\]

Equivalently, $\dim_H(A)$ is the unique value where the Hausdorff measure transitions from infinite to zero.

\section{Regularity of finite-length paths}
\label{sec:regularity}
\label{s: Mink}

\begin{lemma}\label{l:length-below}For every $b>0$ there exists a random constant $C_b$ so that for all paths $\pi$ defined on a time interval $[s,t]\subset[-b,b]$ with $\sup_{r\in[s,t]} |\pi(r)|\le b$, and all $s\le r_1<r_2\le t$ we have
$$
{ \cL(\pi(r_1),r_1;\pi(r_2),r_2)}\ge |\pi_{[r_1,r_2]}|\ge -C_b+|\pi|.
$$
\end{lemma}

\begin{proof}
If $|\pi|=-\infty$, the lower bound is immediate. Otherwise,
$$
|\pi_{[r_1,r_2]}|=
 |\pi|- |\pi_{[s,r_1]}|-|\pi_{[r_2,t]}|\ge |\pi|-{ \cL(\pi(s),s;\pi(r_1),r_1)}-{ \cL(\pi(r_2),r_2;\pi(t),t)}.
$$
{ Zero-duration pieces have length zero and are omitted. By Theorem \ref{t:DL-value-bound} we have that}
$$
{ \cL(\pi(s),s;\pi(r_1),r_1)}\le C_b(r_1-s)^{1/3}{ (\log(e+1/(r_1-s)))^{4/3}}-(\pi(r_1)-\pi(s))^2/(r_1-s)\le C_b'
$$
and a similar bound holds for ${ \cL(\pi(r_2),r_2;\pi(t),t)}$.
\end{proof}

\begin{lemma}\label{l:Holder1/2}
A.s., for every $b>0$ all finite length paths $\pi$ defined on a time interval $[s,t]\subset[-b,b]$ with $\sup_{r\in[s,t]} |\pi(r)|\le b$ have { finite} H\"older-1/2 { seminorm, uniformly when their lengths are bounded below}.
More precisely, with random constants $C_b,C_b'$ depending on $b$ only, { for every $[v,v+u]\subset[s,t]$,}
{ \begin{align}\label{e:Holder1/2-1}
(\pi(v+u)-\pi(v))^2 &\le -u |\pi_{[v,v+u]}| + C_bu^{4/3}(\log(e+1/u))^{4/3} \notag \\
&\le (C_b-|\pi|)u + C_b u^{4/3}(\log(e+1/u))^{4/3}
\le (C_b'-|\pi|)u.
\end{align}}
\end{lemma}
\begin{proof}
By Lemma \ref{l:length-below},
{ \begin{align*}
-C_b+|\pi|&\le |\pi_{[v,v+u]}|
\le \cL(\pi(v),v;\pi(v+u),v+u)\\
&\le -\frac{(\pi(v+u)-\pi(v))^2}{u}
+C_b' u^{1/3}(\log(e+1/u))^{4/3}.
\end{align*}}
where the last bound follows from Theorem \ref{t:DL-value-bound}. Rearranging { and using $0<u\le2b$}, we get the claim.
\end{proof}

\begin{lemma} A.s., for any path $\pi:[s,t]\to\R$ with finite length, and $\alpha>3/2$, we have
$$
\lim_{\eps\to 0}\sum_{i=0}^{N_\eps-1}
|\pi(s+(i+1)\eps)-\pi(s+i\eps)|^\alpha=0,
\qquad N_\eps=\lfloor(t-s)/\eps\rfloor.
$$
\end{lemma}
In fact we could have $\alpha=3/2$ with the right logarithmic correction. 
\begin{proof}
Without loss of generality, we may assume that $\pi$ is on the time interval $[0,1]$. Put $N=\lfloor1/\eps\rfloor$ and let $\Delta_i \pi = |\pi(i\eps )-\pi((i-1)\eps )|$. First let $3/2<\alpha\le2$. By H\"older's inequality, 
\begin{align*}
\br{\sum_{i=1}^{N}\Delta_i \pi^\alpha}^{1/\alpha} \le 
\br{\sum_{i=1}^{N}\Delta_i \pi^2}^{1/2}\eps^{1/2-1/\alpha}
\end{align*}
Summing the landscape bound over the grid intervals and using Lemma \ref{l:length-below} for $\pi_{[0,N\eps]}$ gives
\begin{equation*}
\sum_{i=1}^{N}|\pi(i\eps)-\pi((i-1)\eps)|^2\le (C-|\pi|)\eps+C\eps^{1/3}(\log(e/\eps))^{4/3}.
\end{equation*}
Thus the preceding expression has the upper bound
\[
((C-|\pi|)\eps+C{ (}\log{ (e/}\eps{ ))}^{4/3}\eps^{1/3})^{1/2}
\eps^{1/2-1/\alpha}\to0,\qquad \eps\to0.
\]
 For $\alpha>2$, continuity and the case $\alpha=2$ give
\[
\sum_i|\Delta_i\pi|^\alpha
\le \left(\max_i|\Delta_i\pi|\right)^{\alpha-2}
\sum_i|\Delta_i\pi|^2\longrightarrow0.\qedhere
\]
\end{proof}

\begin{proposition}\label{p: Mink}
Almost surely, for every $b\ge1$ there is a finite random constant $K_b\ge1$, depending only on the landscape and $b$, such that the following holds simultaneously for every road $\pi:[0,1]\to[-b,b]$ and every integer $k\ge2$. Set
\[
m_i=\max_{s\in[i/k,(i+1)/k]}|\pi(i/k)-\pi(s)|,
\qquad h_\pi=(-|\pi|)_+.
\]
Then
\begin{align}
\sum_{i=0}^{k-1}m_i^2
&\le-\frac{|\pi|}{k}+K_b k^{-1/3}(\log(ek))^{4/3},
\label{eq:squared-modulus-explicit}\\
\sum_{i=0}^{k-1}m_i
&\le(K_b+h_\pi)^{1/2}k^{1/3}(\log(ek))^{2/3}.
\label{eq:summed-modulus-explicit}
\end{align}
For every partition $0=t_0<\cdots<t_\ell=1$, with $r_j=t_j-t_{j-1}$, the same constant also satisfies
\begin{equation}\label{eq:partition-energy}
\sum_{j=1}^{\ell}\frac{|\pi(t_j)-\pi(t_{j-1})|^2}{r_j}
\le-|\pi|+K_b\sum_{j=1}^{\ell}r_j^{1/3}(\log(e/r_j))^{4/3}.
\end{equation}
Thus this bound holds simultaneously for all sampled partitions, including partitions chosen from the landscape and the path.
The constants $K_b$ may be chosen nondecreasing in $b$. In particular, for paths with ${\sup_{t\in[0,1]}|\pi(t)|}\le b$ and $|\pi|\ge-h$, the bounds are uniform: replace $h_\pi$ by $h$ in the second estimate, and the first implies
\[
\sum_i m_i^2\le(K_b+h)k^{-1/3}(\log(ek))^{4/3}.
\]
\end{proposition}
\begin{proof}
The landscape bound on $[-b,b]\times[0,1]$ gives a random $A_b$, independent of the choice of path, such that every segment of duration $u>0$ satisfies
\[
\frac{|\pi(t+u)-\pi(t)|^2}{u}
\le-|\pi_{[t,t+u]}|+A_b u^{1/3}(\log(e/u))^{4/3}.
\]
Sum this inequality over the intervals of an arbitrary partition. Additivity of path length gives \eqref{eq:partition-energy}, after enlarging $K_b$ if necessary. To obtain the equal-grid modulus estimates, note that, uniformly for $0<u\le1/k$,
\[
u^{1/3}(\log(e/u))^{4/3}
\le c k^{-1/3}(\log(ek))^{4/3}.
\]
Choose $s_i\in[t_i,t_{i+1}]$ attaining $m_i$, where $t_i=i/k$. Apply the segment estimate on $[t_i,s_i]$ and $[s_i,t_{i+1}]$. If one has zero duration, omit it; if $s_i=t_i$, then $m_i=0$. Adding the two estimates gives
\[
k m_i^2\le-|\pi_{[t_i,t_{i+1}]}|
+cA_b k^{-1/3}(\log(ek))^{4/3}.
\]
Sum over $i$ and use additivity of path length:
\[
k\sum_i m_i^2\le-|\pi|+K_b k^{2/3}(\log(ek))^{4/3},
\]
with $K_b\ge1$ depending only on $A_b$. This proves \eqref{eq:squared-modulus-explicit}. Cauchy--Schwarz now gives
\begin{align*}
\sum_i m_i
&\le\left(k\sum_i m_i^2\right)^{1/2}\\
&\le\left(h_\pi+K_b k^{2/3}(\log(ek))^{4/3}\right)^{1/2}\\
&\le(K_b+h_\pi)^{1/2}k^{1/3}(\log(ek))^{2/3}.
\end{align*}
Choose the landscape constants first on integer-sized boxes, replace them by their running maxima, and use the box of size $\lceil b\rceil$ for general $b$. This makes $K_b$ nondecreasing and the estimates simultaneous for all $b$ and all paths. The uniform bounds under $|\pi|\ge-h$ follow from $h_\pi\le h$ and $k^{-1}\le k^{-1/3}(\log(ek))^{4/3}$.
\end{proof}

We now have all the tools required to prove that every road has Minkowski dimension at most $4/3$ almost surely.
\begin{proof}[Proof of Theorem \ref{t: Mink}]
By shift and scale invariance, it suffices to show this for paths $\pi:[0,1]\to \mathbb R$. By Proposition \ref{p: Mink}, on each time interval $[i/k,(i+1)/k]$, the number of squares of side length $1/k$ that will cover the graph of $\pi$ over that interval is at most ${ \lceil2m_ik\rceil+1}$. So the { total number of boxes needed} to cover the function is bounded by $C' k^{4/3}\log^{2/3}k$. Thus the Minkowski dimension of $\pi$ is at most $4/3$.
\end{proof}

\section{An FKG bound on the environment outside a geodesic}
\label{sec:fkg}

Our next objective is to establish that knowledge of a geodesic in the directed landscape does not increase the environment around it. Analogous results are known in discrete last-passage percolation models (see, for example \cite{BSS19, MSZ25}), though passing to the limit is nontrivial. In the continuum setting, previous works have often relied on pre-limiting inequalities, taking limits only in specific applications. We introduce a new formulation that holds in the directed landscape. This formulation is { equivalent to} Theorem \ref{t:FKG-intro}.

\begin{theorem}\label{t:FKG}
Let $u\in\dR$, let $\gamma$ be the $u$-geodesic in the directed landscape $\cL$. Let \(\Pi_\gamma(v)\) be the set of \(v\)-paths with graphs disjoint from that of \(\gamma\), and define
\[
\mathcal{L}_\gamma(v) := 
\begin{cases}
\sup_{\pi \in \Pi_\gamma(v)} |\pi| & \text{if } \Pi_\gamma(v) \neq \emptyset, \\
-\infty & \text{else}.
\end{cases}
\]
There exists a coupling of $\cL$ and another directed landscape $\cL'$ so that $\cL_\gamma \le \cL'$ and $\gamma$ is independent of $\cL'.$
\end{theorem}

A direct consequence is the classical formulation.

\begin{corollary}
 For any bounded measurable function $f$, non-decreasing in the second coordinate, $\EE\bs{f(\gamma, \cL_{\gamma})}\le \EE \bs{f(\gamma, \cL')}$.
\end{corollary}

However, it seems that the formulation of Theorem \ref{t:FKG} behaves much better under taking limits. The first step of the proof is to show the analogous theorem for exponential last passage. 

\subsection{Exponential last passage percolation}
\label{subsec:expLPP}

We record the basic definitions of exponential last passage percolation on the integer lattice; this discrete model will serve as the setting for the coupling argument that yields the FKG domination in the continuum. Let \(\omega=(\omega_{i,j})_{i,j\in\bZ}\) be a field of i.i.d.\ exponential random variables with rate one. Define the set of \textbf{admissible endpoint pairs} by
\[
\mathcal{P} := \bigl\{ \bigl( (x,y),(x',y') \bigr) \in \bZ^2 \times \bZ^2 : x\le x',\ y\le y' \bigr\}.
\]
For any \((p,q)\in\mathcal{P}\), a \emph{directed path} \(\pi\) from \(p=(x,y)\) to \(q=(x',y')\) is a finite sequence of lattice points \((i_k,j_k)_{k=0}^m\) satisfying \((i_0,j_0)=p\), \((i_m,j_m)=q\), and each consecutive pair differing by either a unit right step or a unit up step. Write \(\Pi(p,q)\) for the finite set of all such paths.

The \textbf{length} of a path \(\pi\) is the sum of the weights along its vertices:
\[
|\pi| := \sum_{(i,j)\in\pi} \omega_{i,j}.
\]
The \textbf{last passage distance} is then
\[
d(p,q) := \max_{\pi\in\Pi(p,q)} |\pi|.
\]
Any path attaining the maximum is called a \textbf{geodesic}.

Fix an arbitrary directed path \(\eta\) from \(p\) to \(q\) with \((p,q)\in\mathcal{P}\). For another admissible pair \(v=(r,s)\in\mathcal{P}\), let \(\Pi_\eta(v)\) be the set of all directed paths from the starting point of \(v\) to the ending point of \(v\) whose vertex set is disjoint from the vertex set of \(\eta\). Define the \textbf{restricted last passage distance}
\begin{equation}\label{e: restricted-lp}
d_\eta(v) :=
\begin{cases}
\max_{\pi\in\Pi_\eta(v)} |\pi| & \text{if }\Pi_\eta(v)\ne\emptyset,\\
-\infty & \text{otherwise}.
\end{cases}
\end{equation}

\begin{theorem}
\label{t:FKG-exp}
Let \((p,q)\in\mathcal{P}\). Let \(\gamma\) be a geodesic from \(p\) to \(q\). There exists a coupling of the last-passage distance field \(d\) and a copy \(d'\ed d\) such that
\[
d_\gamma \le d' \qquad\text{a.s.}
\]
and \(\gamma\) is independent of \(d'\).
\end{theorem}

\begin{proof}[Proof of Theorem \ref{t:FKG-exp}]
Let \((p,q)\in\mathcal{P}\) and let \(\Pi(p,q)\) be the finite set of all directed paths from \(p\) to \(q\). Consider the partial order on the product space
\[
\Pi(p,q) \times { (\mathbb R\cup\{-\infty\})^{\mathcal P}}
\]
defined by \((\pi_1,d_1) \preceq (\pi_2,d_2)\) if and only if \(\pi_1=\pi_2\) and \(d_1\le d_2\).

Let \({ f:\Pi(p,q)\times(\mathbb R\cup\{-\infty\})^{\mathcal P}\to\R}\) be bounded, measurable, and nondecreasing with respect to this partial order. Let \(d\) be the last-passage distance field on \(\mathcal{P}\) (induced by i.i.d.\ exponential weights), let \(\gamma\) be the geodesic from \(p\) to \(q\), and let \(d_\gamma\) be the associated restricted last-passage distance as in \eqref{e: restricted-lp}. Let \(d'\) be an independent copy of \(d\). For any fixed path \(\pi\in\Pi(p,q)\), write \(d_\pi\) (resp.\ \(d'_\pi\)) for the restricted last-passage distance with respect to the path \(\pi\) computed from \(d\) (resp.\ from the independent copy \(d'\)).

Fix any path \(\pi\in\Pi(p,q)\). On the event \(\{\gamma=\pi\}\) we have \(d_\gamma = d_\pi\). 

We have
\[
\EE\bigl[(-f(\pi,d_\pi))\,1_{\{\gamma=\pi\}}\bigr] \ge \EE\bigl[(-f(\pi,d_\pi))\bigr]\,\EE\bigl[1_{\{\gamma=\pi\}}\bigr].
\]
by the FKG inequality, since  each factor in the first expectation is nonincreasing in the weights off $\pi$ and the negative of the weights on $\pi$.

The right-hand side equals \(\EE\bigl[(-f(\pi,d'_\pi))\,1_{\{\gamma=\pi\}}\bigr]\), because \(d'\) is an independent copy of \(d\) (hence \(d'_\pi\) has the same law as \(d_\pi\)) and the indicator depends only on \(d\). Rearranging gives
\[
\EE\bigl[f(\pi,d_\pi)\,1_{\{\gamma=\pi\}}\bigr] \le \EE\bigl[f(\pi,d'_\pi)\,1_{\{\gamma=\pi\}}\bigr].
\]
Summing over the finitely many paths \(\pi\in\Pi(p,q)\) yields
\[
\EE \bs{f(\gamma,d_\gamma)} \le \EE \bs{f(\gamma,d'_\gamma)} \le \EE \bs{f(\gamma,d')},
\]
where the last step uses the elementary inequality \(d'_\gamma\le d'\) (any path that avoids \(\gamma\) is a valid path in the full environment, so its maximum is at most the unrestricted last-passage distance \(d'\)).

The inequality \(\EE [f(\gamma,d_\gamma)]\le\EE[ f(\gamma,d')]\) then holds for all bounded non-decreasing test functions. By Strassen's theorem (Theorem \ref{t:strassen}) there exists a coupling \(C_1\) of the laws of \((\gamma,d_\gamma)\) and \((\gamma',d')\) (with \(\gamma'\) independent of \(d'\)) such that \((\gamma,d_\gamma)\preceq(\gamma',d')\) almost surely. In particular, by definition of the partial order, \(\gamma=\gamma'\) a.s.\ in this coupling.

The pair \((\gamma,d_\gamma)\) is a deterministic measurable function of the field \(d\). This gives a trivial coupling \(C_2\) of \(d\) with \((\gamma,d_\gamma)\). The gluing lemma (Lemma \ref{l:gluing}) allows us to realize \(C_1\) and \(C_2\) on a common probability space. The resulting coupling of \(d\) and \(d'\) satisfies \(d_\gamma\le d'\) almost surely and \(\gamma\) is independent of \(d'\).
\end{proof}

We need two simple lemmas about the directed landscape. 
\begin{lemma}[Splitting for uniqueness]\label{l:geodesic-split}
 For every $u=(x_1,t_1;x_4,t_4)\in \dR$ and $u$-geodesic $\gamma$ { there exist $t_2,t_3$} with $t_1<t_2<t_3< t_4$ so that for $i=1,2,3$ the $u_i=(\tilde \gamma(t_i),\tilde\gamma(t_{i+1}))$-geodesic is unique. 
 \end{lemma}
\begin{proof}
By Theorem \ref{t:finite-geodesics} there are at most finitely many $u$-geodesics for every $u$. Moreover,  the no bubbles property Theorem \ref{t:no-bubbles} says that for every two $u$-geodesics, their set of intersection times is of the form $\{t_1,t_4\}\cup I$ for some closed interval $I$. Thus, the $u_2$-geodesic is always unique, or else we could create a different intersection pattern by swapping the $t_2-t_3$ section of $\gamma$ to the alternative. 

As $t_2\downarrow t_1$, then the  $u_1$-geodesic must eventually be unique, or else we could create  infinitely many $u$-geodesics by swapping the $t_1-t_2$ section of $\gamma$ to the alternative. 

By a symmetric argument, the same holds for $u_3$ as $t_3\uparrow t_4$. 
\end{proof}

\begin{lemma}[Piecewise geodesic approximation]\label{l:PGA}
The following holds almost surely in $\cL$. Let $u\in \dR$ and  $\pi$ be a finite length $u$-path. Then there { exist $u$-paths} $\pi_n$ so that
\begin{enumerate}
    \item For each $n$, $\pi_n$ is a concatenation of $v$-geodesic segments for $v\in V_n\subset \dR$, a finite subset.
    \item $\pi_n\to\pi$ uniformly
    \item $\lim_{n\to \infty}|\pi_n|=|\pi|$. 
    \item For each $v\in V_n$ there is a unique $v$-geodesic in $\cL$.
\end{enumerate}
\end{lemma}
\begin{proof}
Recall the  definition \eqref{e:length} of path length
\begin{equation}\label{e:length2}
|\pi|=\inf\sum_{i=1}^k\mathcal{L}(\tilde \pi(t_{i-1}),\tilde \pi(t_i)),
\end{equation}
where the infimum runs over all finite subdivisions \(s=t_0<\dots<t_k=t\). The triangle inequality implies that the sum does not increase if we replace a subdivision by a refinement. Thus we may take a sequence of more and more refined subdivisions with mesh going to zero along which the sum converges. 

Let $\pi_n$ be the path that agrees with $\pi$ at the times of the $n$th subdivision but consists of geodesics in between. { Then $|\pi_n|\to |\pi|$. The fixed endpoints and the lower bound on length keep these paths in a common compact region by Theorem \ref{t:DL-value-bound}. Since paths} with length bounded below form a { uniformly} equicontinuous family (see Lemma \ref{l:Holder1/2}), $\pi_n$ must have subsequential limits in the uniform norm, and these all agree with $\pi$ at all subdivision times, and hence everywhere. This construction achieves $1-3$. 

Finally every $v=(x,s;y,t)$-geodesic can be split in 3 parts at intermediate points so that $4$ also holds by Lemma \ref{l:geodesic-split}. 
\end{proof}

\subsection{Convergence to the landscape and proof of Theorem \ref{t:FKG}}

\cite{dauvergne2021scaling} construct a coupling of the directed landscape \(\mathcal{L}\) with a sequence of rescaled exponential last-passage percolation models \(d_n\) (extended to upper semi-continuous functions on \(\dR\)) such that the following hold almost surely:
\begin{enumerate}[label = (\arabic*)]
\item \(d_n\to\mathcal{L}\) in the graph topology (this implies uniform convergence on compact subsets of \(\dR\)).
\item In this setting, Theorem 8.7 in \cite{dauvergne2021scaling} combined with (1) implies that for every \(u\in\dR\) whose \(u\)-geodesic in \(\mathcal{L}\) is unique, the rescaled geodesic of \(d_n\) from lattice approximations to these endpoints converges uniformly (as a function of time) to the \(u\)-geodesic of \(\mathcal{L}\).
\end{enumerate}
For the definition of graph convergence and the exact rescaling parameters, we refer to \cite{dauvergne2021scaling}.

This joint convergence is the key tool that allows discrete FKG-type domination statements to pass to the continuum limit.

\begin{proof}[Proof of Theorem \ref{t:FKG}]
Fix $u\in \dR$. Since $u$ is deterministic, its $\cL$-geodesic is almost surely unique
by \cite[Theorem 12.1]{DOV}; denote it by $\gamma$. 
We start with the rescaled version $(d_n,d_n')$ of the coupling in Theorem \ref{t:FKG-exp}, where $\gamma_n$ is a geodesic in $d_n$ with endpoints approximating $u$, $(d_n)_{\gamma_n}\le d_n'$, and $d_n'$ is independent of $\gamma_n$. 

The convergence theorem implies that $(d_n,d_n')$ is tight in the product graph topology. Passing to a subsequence and using Skorokhod embedding, we may assume that
\[
(d_n,d_n')\longrightarrow(\mathcal L,\mathcal L')
\]
almost surely, where both marginals are directed landscapes. We need to prove that $\cL'$ is independent of $\gamma$ and that
\begin{equation}\label{e: LgLs}
\mathcal L_\gamma(w)\le\mathcal L'(w),\qquad w\in\dR.
\end{equation}

The geodesic convergence gives $\gamma_n\to\gamma$ uniformly, and hence $(d_n',\gamma_n)\to(\mathcal L',\gamma)$. Since independence is preserved in the limit, we have shown that $\mathcal L'$ and $\gamma$ are independent. 

Now let $w\in \dR$, $\varepsilon> 0$. If $\cL_\gamma(w)>-\infty$, then there exists a path $\pi$ disjoint from $\gamma$ so that $|\pi|_{\cL}\ge \cL_\gamma(w)-\eps$. Moreover, by Lemma \ref{l:PGA}, we may assume that $\pi$ is a concatenation of  finitely many $\mathcal L$-geodesic segments $\pi_v$ with endpoint pairs $v\in V$, and each $\pi_v$ is the unique $v$-geodesic.  We choose lattice approximations $v_n$ to the pairs $v\in V$, using the same lattice point at each common endpoint. By Theorem 8.7 in \cite{dauvergne2021scaling}, $d_n(v_n)\to\cL(v)$ and the $v_n$-geodesic $\pi_{n,v}$ in $d_n$ converges to $\pi_v$ uniformly. Concatenating these segments gives a path $\pi_n$ with endpoint pair $w_n\to w$ and $|\pi_n|_{d_n}\to|\pi|_\cL$. 
The graphs of $\pi$ and $\gamma$ are disjoint compact sets, so they
have positive separation. Since $\pi_n\to\pi$ and $\gamma_n\to\gamma$
uniformly, $\pi_n$ is disjoint from $\gamma_n$ for all 
large $n$. Consequently,
$|\pi_n|_{d_n}\le (d_n)_{\gamma_n}(w_n)\le d_n'(w_n)$. By graph convergence,
\[
\cL'(w)=\lim_{n\to\infty}d_n'(w_n)
\ge\lim_{n\to\infty}|\pi_n|_{d_n}
=|\pi|_\cL\ge\cL_\gamma(w)-\eps,
\]
showing \eqref{e: LgLs} as we let $\eps\downarrow 0$. This completes the proof. \end{proof}

\begin{corollary}\label{c:FKG}
In the above coupling, we have 
$\cL_\gamma \le \cL'_\gamma.$ 
For any $\cL$-path  $\pi$ defined on $[s',t']$  that intersects $\gamma$ first at time $s$ and last at time $t$ we have
$$
|\pi|_\cL\le \cL(\tilde\pi(s),\tilde\pi(t))+|\pi_{[s',s]}|_{\cL'}+|\pi_{[t,t']}|_{\cL'}\le  \cL(\tilde\pi(s),\tilde\pi(t))+\cL'(\tilde\pi(s'),\tilde\pi(s))+\cL'(\tilde\pi(t),\tilde\pi(t'))
$$
\end{corollary}
\begin{proof}
The first claim follows from the definition of $\cL_\gamma$ and $\cL'_{\gamma}$ and the fact that $\cL_\gamma \le \cL'$.
For the second claim, we prove the case $s'<s$, $t<t'$, as the other cases are similar. 
By the definition of length, 
$$
|\pi|_{\cL}\le \cL(\tilde\pi(s-\eps),\tilde\pi(t+\eps))+|\pi_{[s',s-\eps]}|_{\cL}+|\pi_{[t+\eps,t']}|_{\cL}
$$
since the two path segments on the right hand side are disjoint from $\gamma$, we can upper bound $|\cdot|_{\cL}$ by $|\cdot|_{\cL'}$. Now take a limit as $\eps \to 0$. To get the first inequality, we use that  $\cL$ is continuous away from $s=t$, upper semicontinuous even at $s=t$, and that $|\pi_{[s',s-\eps]}|_{\cL'}\to |\pi_{[s',s]}|_{\cL'}$, which follows from the definition of path length. The second inequality follows directly from the definition of path length. 
\end{proof}

\section{Good detours are unlikely}
\label{sec:detour}

The goal of this section is to show that detours of good length from a geodesic are unlikely. We start with some preliminary estimates in Section \ref{ss:detour-prelim}. We treat nearby detours in Section \ref{ss:detour-nearby}, far away detours in Section \ref{ss:detour-far}, and put everything together in Theorem \ref{t:detour} of Section \ref{ss:detour-main}.

\subsection{Preliminary estimates}\label{ss:detour-prelim}

We start with a lemma about the directed landscape that will be used in the next proposition. Let $[a\pm \chi]$ denote the interval $[a-\chi,a+\chi]$, and we extend $\cL$ to sets by taking supremum of $\cL$ distance over points in the set.

\begin{lemma}\label{l:LD}
There exists $\eps>0$, $m>0$ and $c,d>0$ so that for all $\sigma>0$ and $n\ge 1$ we have
$$
\PP(X_1+\ldots+X_n\ge -\sigma mn)\le c e^{-dn}$$ whenever the $X_{i}$ are independent with law 
$
\mathcal 
L([a_i\pm \eps_i],t_i;[b_i\pm \delta_i],t_i+\sigma^3)$, and $\eps_i,\delta_i\le \sigma^2\eps$.
\end{lemma}
\begin{proof} By KPZ scaling and shift-invariance, we may assume $\sigma=1$, $t_i=b_i=0$. By increasing all intervals to length $2\eps$, we may assume  $\eps_i=\delta_i=\eps$. First consider the case  $a_i=0$.
By the continuity of $\cL$, as $\eps\downarrow 0$, 
$$
H_\eps=\cL([\pm \eps],0;[\pm \eps],1)\to \cL(0,0;0,1),
$$
and decreases as $\eps\downarrow0$.
These variables are uniformly integrable: in fact they  have uniform exponential tail bounds by Theorem \ref{t:DL-value-bound}.  Since the GUE Tracy-Widom distribution  of $\cL(0,0;0,1)$ has negative expectation (\cite{tracy1994level}), we can fix $\eps>0$ so that $m:=-\EE [H_{5\eps}/2]>0$. Theorem \ref{t:DL-value-bound} implies that the moment generating function for $H_{5\epsilon}$ exists on $\mathbb R$, and therefore the claim follows by Cramer's large deviation theorem.

For the general case, it suffices to show the stochastic domination  
$
X_i\preceq H_{5\eps}.
$
By shear  invariance  as processes in $x,y$
\begin{equation}\label{e:shift}
 \cL(x+a,0;y,1) \ed  \cL(x,0;y,1) -a(a+2(x-y))
\end{equation}
so when $4\eps\le |a|$, the second term is negative for all $x,y\in [\pm \epsilon]$, and we have the stochastic domination 
$
\cL([a\pm \eps],0;[\pm \eps],1) \preceq H_{\eps}\le H_{5\eps}.
$
When $4\eps \ge |a|$, by monotonicity
$\cL([a\pm \eps],0;[\pm \eps],1) \le  H_{5\eps}$.
\end{proof}

The next proposition controls near maxima of Brownian motion with drift on a finite interval.

\begin{proposition}[Near-maxima of Brownian motion with drift]\label{p:BM-max}
Let $B$ be Brownian motion with variance $\sigma^2$ and drift $\nu$ defined on an interval $[x\pm2\rho]$. Let $M,Z$ denote the  maximum and its location (first if there are ties) of $B$. Then for any $\delta>0$ we have
\begin{align}\label{e:MZ1}
&\PP\Big(|Z-x|\le \rho\Big) \le 2\exp\left(-\frac{\nu^2 \rho}{3\sigma^2}\right),
\\ \label{e:MZ3}
&\PP\Big(\sup_{[x\pm \rho]\setminus [Z\pm y]} B\ge M-\delta\sigma\big\vert \,|Z-x|\le \rho\Big)\le 12\delta/\sqrt{y}.
\end{align}
\end{proposition}
\begin{proof}
By Brownian scaling and translation invariance, we may assume $\sigma=1$, $x=0$, $\rho=1$. For \eqref{e:MZ1}, by symmetry it suffices to consider the case $\nu> 0$. Assuming this, after extending $B$ to the entire negative real line, we have that
$$
S=\sup_{(-\infty,1]}B-B(1),
\qquad S'=B(2)-B(1)$$
are independent,  and $|Z|\le 1$ implies $S'\le  S$. 
With standard exponential and normal variables $T$, $N$ we have $S\sim T/(2\nu), S'\sim \nu+N$, so for any $a\in [0,1]$
\begin{align*}
\PP(|Z|\le 1) &\le \PP(T/(2\nu)\ge N+\nu) \\&\le \PP(T/(2\nu)\ge \nu a)+\PP(-N\ge \nu(1-a))\le e^{-2a\nu^2}+e^{-\nu^2(1-a)^2/2},
\end{align*} 
giving \eqref{e:MZ1} after optimizing for $a$. 

For \eqref{e:MZ3}, first assume $\nu=0$. By the Williams decomposition, conditionally on $Z$, the law of $(M-B(t), t\in [Z,2])$ given $Z$ is a Brownian meander on $[Z,2]$.  Similarly the conditional law of $(M-B(t),t\in[-2,Z])$
is a Brownian meander run backwards in time on $[-2,Z]$, and the two meanders are independent. 

Now we move to drifted, general $\nu$ case, then by Girsanov's theorem, the Radon-Nikodym derivative of the conditional law is proportional to 
$\exp(\nu(B(2)-B(-2)))$, which factors as $\exp(\nu(B(2)-M))\exp(\nu(M-B(-2)))$. In particular, the conditional laws given $Z$ are still independent meanders, just biased by a function of their endpoints. (The biasing also changes the marginal law of $Z$, but this is not relevant to our argument.)

Next, we focus on the biased meander $M-B$ on $[Z,2]$. 
Since a meander is just a Bessel-3 bridge given its endpoints, $M-B$ on $[Z,2]$  stochastically dominates $e$, a Brownian excursion on $[Z,2]$.

Next, we restrict our attention to the subinterval $[Z,1]\subset [Z,2]$. Given $Z$ and $e(1)$, the law of $e$ on $[Z,1]$ is that of a Bessel bridge, and given $Z$, $e(1)$ has the same law as $\alpha\chi_3$, where $\chi_3$ has Chi-distribution with $3$ degrees of freedom, and $\alpha=\sqrt{(1-Z)/(2-Z)}$. 

Next, we compare the excursion $e$ on $[Z,1]$ to $R$ restricted to $[Z,1]$, where $R$ is a standard 3-dimensional Bessel process defined on $[Z,\infty)$.  Conditionally on their values at $1$ they are both Bessel-3 bridges, so it suffices to compare the law of their endpoints, $e(1)\sim \alpha \chi_3$ and $R(1)\sim \beta\chi_3$ with  $\beta=\sqrt{1-Z}$.
The Radon-Nikodym derivative of $\alpha \chi_3$ with respect to  $\beta\chi_3$ at $x$ is given by 
$$
(\beta/\alpha)^3\exp(-x^2(\alpha^{-2}-\beta^{-2})/2)\le (\beta/\alpha)^3\le 3^{3/2}.
$$
on $|Z|\le 1$, since then we have $1/\sqrt{3}\le\alpha/\beta \le 1$. 
So on $\{|Z|\le 1\}$ we have
\begin{align*}
\PP&\Big(\sup_{[Z+y,1]} B\ge M-\delta\;|\,Z\Big)
\le \PP\Big(\inf_{[Z+y,1]}e\le \delta\;|\;Z\Big)\\
&\le 3^{3/2}\PP\Big(\inf_{[Z+y,1]}R\le\delta \;|\;Z\Big)
\le 3^{3/2}\PP\Big(\inf_{[Z+y,\infty)}R\le \delta\;|\;Z\Big).
\end{align*}
The last probability depends on $\delta$ and $y$, but not on $Z$. Since the  law of the future minimum of $R$ given $R(t)$ is uniform on $[0,R(t)]$, the last probability equals
\begin{align*}
\PP\Big(\inf_{[Z+y,\infty)}R\le \delta\Big)&=E[1\wedge \delta/R(Z+y)]=E[1\wedge (\delta/(\sqrt{y}\chi_3))]\\&\le \delta/{\sqrt{y}}\,%
E[1/\chi_3]={\sqrt{2/\pi}\,\delta/\sqrt y}.
\end{align*}
A symmetric argument on the left side of $Z$ gives that on $\{|Z|\le 1\}$,
$$\PP\Big(\sup_{[-1,Z-y]} B\ge M-\delta\;|\;Z\Big)\le 2\times3^{3/2}/\sqrt{\pi}\times \delta/\sqrt{y},$$
and the two bounds imply \eqref{e:MZ3}.
\end{proof}

The next lemma controls the hazard rate of standard normal random variables.

\begin{lemma}[Hazard rate bound for  normals]\label{l:normal-delta}
For  $Z$ normal with mean $b$ and variance $1$ and any $a,z\in \R$ and $\delta>0$, 
\begin{equation}\label{e:Gauss}
    \PP(Z\in [a,a+\delta]|Z>z)\le c((z-b)\vee 1)\delta\,.
\end{equation}
\end{lemma}
\begin{proof}
Note that, by subtracting the mean from $Z$, we can without loss of generality assume that $b=0$. We have also assumed that $\delta<1/2$ . If $z\le 2$, the probability is bounded by \[\PP(Z\in [a,a+\delta])\PP(Z>2)^{-1}\le c\delta.\] 
Next, assume that $z>2$. Comparing with $a'=z$ on the left shows that we may also assume $a\ge z$. Then by standard bounds on Gaussian distribution functions (see, for example, \cite{morters2010brownian} Lemma $12.9$)
\[\PP(Z>z)\ge (2z)^{-1}\phi(z)\,,\]
where $\phi$ is the density function of standard Gaussian. Then since $\phi(z)$ is decreasing for $z>0$, we have 
\[
\frac{\PP(Z\in [a,a+\delta])}{\PP(Z>z)}\le 2z\phi(z)^{-1}{\PP(Z\in [a,a+\delta]})\le 2z\phi(z)^{-1}\phi(a)\delta\le 2z\delta. \qedhere
\]
\end{proof}

\subsection{Good nearby detours}\label{ss:detour-nearby}

The following proposition shows that nearby detours are unlikely.

\begin{proposition}[Good nearby detours are unlikely]\label{p:detour-fkg}
Let $\gamma$ be a geodesic from $(x,r')$ to $(y,t')$. Let $(r,r+\sigma^3)\subset (r',t')$. Let $\Pi$ be the set of paths $\pi$ with the same endpoints that agree with $\gamma$ on $(r',t')\setminus (r-2\sigma^3,r+3\sigma^3)$, are disjoint from $\gamma$ on 
$(r,r+\sigma^3)$, and  satisfy
{
\begin{equation}\label{e:detour-fkg}\max_{i=1,\dots,n-1}|\gamma(t_i)-\pi(t_i)|\le \eps_0 \sigma^2/n^{2/3}, \qquad  t_i=r+i\sigma^3/n.
\end{equation}
}
Then
\begin{equation}
\PP\Big(\inf_{\pi \in \Pi}\Delta(\pi)  \le c_1\sigma n^{2/3}\big)\le c_2e^{-dn}
\end{equation}
with universal constants $c_1,c_2,\eps_0,d>0$ (where $\eps_0$ is the constant from Lemma \ref{l:LD}).
\end{proposition}

\begin{proof}{ By shear and scaling, we may assume $r=0$, $\sigma=1$, $x=y=0$. 
By enlarging $c_2$, the cases $n<4$ follow from the bound $\PP(\cdot)\le1$, so we assume $n\ge4$. Put $a=t_1$ and $b=t_{n-1}$. We write
\begin{align}\label{e:4terms}
\begin{split}
\Delta(\pi)=|\gamma|-|\pi|&=|\gamma_{[a,b]}|\\&
-|\pi_{[a,b]}|\\
&+
|\gamma_{[r',a]}|-|\pi_{[r',a]}|\\
&+|\gamma_{[b,t']}|-|\pi_{[b,t']}|
\end{split}
\end{align}
Denote the  terms in the expression above so that it reads $Y(\pi)=Y_1+Y_2(\pi)+Y_3(\pi)+Y_4(\pi)$. 
Let $Y=\inf_{\Pi}Y(\pi)$, and define $Y_i$ analogously. 
We will use the union bound $$
\PP(Y\le c_1n^{2/3})\le \PP(Y_1\le -c_1n^{2/3})+\PP(Y_2\le 4c_1n^{2/3})+\PP(Y_3\le -c_1n^{2/3})+\PP(Y_4\le -c_1n^{2/3}).
$$
We have $\PP(Y_1\le-c_1n^{2/3})\le ce^{-dn}$ by Theorem \ref{t:geod-thm-DSV}.
For the second term, by the definition of path length, for any path $\pi$
$$
|\pi_{[a,b]}| \le \sum_{i=2}^{n-1} \cL(\tilde\pi( t_{i-1});\tilde \pi(t_i))\,.
$$
where $\tilde\pi(t)=(\pi(t),t)$.
Let $\cL_\gamma(u)=\sup |g|$ over all paths $g$ with endpoints $u$ having no point in common with $\gamma$.  Recall that $\cL_{\gamma}$ .to sets means the supremum, as is usual with distances of negative sign. Let $\chi=\eps_0/ n^{2/3}$ for the $\eps_0$  of Lemma \ref{l:LD}.
Since every $\pi\in \Pi$ is $\chi$-close to $\gamma$ at times $t_i$, we have the upper bound
\begin{align*}-Y_2&\le
\sum_{i=2}^{n-1} \cL_{\gamma}([\gamma(t_{i-1})\pm \chi],t_{i-1};[\gamma(t_{i})\pm \chi],t_{i})
\le
\sum_{i=2}^{n-1} \cL'([\gamma(t_{i-1})\pm \chi],t_{i-1};[\gamma(t_{i})\pm \chi],t_{i})
\end{align*}
In the last inequality, we used Theorem \ref{t:FKG}, which provides  $\cL'\ge \cL_\gamma$, a directed landscape independent of $\gamma$. Write
\[
S:=\sum_{i=2}^{n-1}
\cL'([\gamma(t_{i-1})\pm\chi],t_{i-1};
       [\gamma(t_i)\pm\chi],t_i).
\]
Thus $-Y_2\le S$. Conditional on $\gamma$, the summands in $S$
are independent.

We apply Lemma \ref{l:LD} with $\sigma=n^{-1/3}$. Choosing $c_1=m/8$, we have $4c_1n^{2/3}\le\sigma(n-2)m$ and
\begin{align*}
\PP(Y_2\le 4c_1n^{2/3}\mid\gamma)
&\le\PP(-Y_2\ge-\sigma m(n-2)\mid\gamma)\\
&\le \PP(S\ge-\sigma m(n-2)\mid\gamma)
\le ce^{-dn}.
\end{align*}
Since the bound does not depend on $\gamma$, it holds  unconditionally. 
The last two terms in $\eqref{e:4terms}$ are symmetric, so it { suffices to bound} the first of the two. 
Set $\tau=(-2)\vee r'$. Since for $\pi\in \Pi$ we have $\gamma(\tau)=\pi(\tau)$, we have
$$
-Y_3\le  \cL(\tilde\gamma(\tau);[\gamma(a)\pm \chi],a)-|\gamma_{[\tau,a]}|
$$

Let $\gamma'$ be a geodesic from $\gamma(\tau)$ to $[\gamma(a)\pm \chi]\times\{a\}$, and let $T$ be the last intersection time of $\gamma$ and $\gamma'$. Then we can write the last display as 
$$
-Y_3\le |\gamma_{[\tau,T]}|+|\gamma'_{[T,a]}|-|\gamma_{[\tau,a]}|=|\gamma'_{[T,a]}|-|\gamma_{[T,a]}|
$$
in the last equality, we removed the duplicated segments. The varying part is contained in $[\tau,a]$, whose length is at most $3$. Next, we bound
$$
|\gamma'_{[T,a]}|\le \cL'(\tilde \gamma(T),\tilde \gamma'(a))\le 
\cL'(\tilde \gamma(T);[\gamma(a)\pm \chi],a)
= \cL''(\tilde \gamma(T)-\tilde \gamma(a);[\pm\chi],0)
$$
To apply domination when $\gamma'$ meets $\gamma$ at $T$, we first restrict to $[T+u,a-u]$ and let $u\downarrow0$. Continuity of $\cL'$ and of restricted length for finite-length paths gives the first inequality. If $T=a$, the segment has length zero and needs no estimate. The domination is provided by Corollary \ref{c:FKG}, where $\cL'$ is a directed landscape independent of $\gamma$, and $\cL''$ is a shifted version (still independent of $\gamma$).
By the landscape bound, Theorem \ref{t:DL-value-bound} we have 
$$
 \cL''(\tilde \gamma(t)-\tilde \gamma(a);[\pm\chi],0) \le \sup_{x\in \R, t\in [-3,0)} \cL''(x,t;[\pm\chi],0)\le C
$$
with $\PP(C>a)\le ce^{-da^{3/2}}$. 
Together with the bound on the tails of 
$$-|\gamma_{[T,a]}|\le -\inf_{t\in[\tau,a]} \cL(\tilde \gamma(t),\tilde \gamma(a))
$$
from Corollary \ref{c:local-modulus}, by scaling and using $a-\tau\le3$, we get
\[
\PP(Y_3\le -c_1 n^{2/3})\le ce^{-dn}.
\qedhere\]
}
\end{proof}

\subsection{Good far away detours}
\label{ss:detour-far}

In this section, we will control far-away detours via Proposition \ref{p:detour-far}. The main ingredients are the near-maximum bound for Brownian motion in Proposition \ref{p:BM-max}, and a version of the Brownian Gibbs property introduced by \cite{hammond2022brownian}. For this version, we need the notion of affine shifted bridge version of functions, and a fact about Brownian bridges. 

{ For a real-valued} function $f$ on an interval $[a,b]$, let us define the {\bf affinely shifted bridge version} that is zero at both endpoints, $f^{[a,b]}:[a,b]\to \R$ by
\begin{equation}\label{e: affine shift}
f^{[a,b]}(x) := f(x) - \frac{x-a}{b-a}\cdot f(b) - \frac{b-x}{b-a}\cdot f(a)
\end{equation}

The following elementary lemma is true for Brownian bridges.
\begin{lemma}\label{l:BB}
Let $B$ be a Brownian bridge of variance 1 on $[-2,2]$ with $B(-2)=a$, $B(2)=b$.  Then $B(0)$ is a Gaussian random variable with mean $\frac{a+b}{2}$ and variance $1$. Moreover, if $B_1=B\big|_{[-2,0]}$ and $B_2=B\big|_{[0,2]}$ then $B_1^{[-2,0]}, B_2^{[0,2]}$ and $B(0)$ are independent. 
\end{lemma}
The next proposition shows that far-away detours are unlikely.

\begin{proposition}[Good far away detours are unlikely]\label{p:detour-far}
Let $\gamma$ be a geodesic from $(x,r')$ to $(y,s')$. Let $c_1\ge 1$,  let
$\sigma \le c_1\min(r-r',s'-r)^{1/3}$ let $\delta>0$, and let $A$ be the event that there exists a path $\pi$ with the same endpoints
with 
\begin{equation}\label{e:detour-far}|\gamma(r)-\pi(r)|\ge \sigma^2 ,\qquad  \mbox{and} \qquad  \Delta(\pi) \le \delta\sigma.
\end{equation}
Then there exists $c_2$ depending on $c_1$ only so that
$
\PP(A)\le c_2\delta \log(1/\delta+1).
$
\end{proposition}
\begin{proof}
By shift, shear and scale invariance and time-reversal symmetry, we can assume $r=0$, $x=y=0$, $r'=-1$ and $s'\ge 1$. This scaling ensures that $Z:=\gamma(r)=\gamma(0)$ is typically of order 1. Note also that writing $\sigma'=\sigma/c_1$, $\delta'=c_1\delta$ the $c_1=1$ case implies the general $c_1$ case, so we may assume $c_1=1$, and $\sigma\le 1$. Finally, by adjusting $c_2$ we may also assume $\delta\le 1/2$.

Let $\cA(\cdot)=\cL(0,-1;\cdot,0)$ and $\cA'(\cdot)=\cL(\cdot,0;0,s')$, and let $\mathcal B=\cA+\cA'$. Then $Z=\arg\max \cB$. Let
\[
X\in\operatorname*{arg\,max}_{|z-Z|\ge\sigma^2}\cB(z),
\]
choosing the smallest maximizer if there are ties.

The first case we consider is $|X-Z|\le 3$. When this and $A$ holds,  and $ Z \in[i-1,i+1]$, we also have the following:
$$
\cB(Z)= \sup_{[i-8,i+8]}\cB =\sup_{[i-4,i+4]}  \cB \le \sigma\delta+\sup_{[i-4,i+4]\setminus [Z\pm \sigma^2]}\cB=\sigma\delta + \cB(X).
$$
Now $\cB-\cB(i)$ on the interval $[i-8,i+8]$ is absolutely continuous with respect to two-sided Brownian motion (anchored at $i$) with %
drift { $\nu_i=-2i-2i/s'$, where $2|i|\le|\nu_i|\le4|i|$} with Radon-Nikodym derivative bounded by a universal constant, see \cite{dauvergne2024wiener}.

By Proposition \ref{p:BM-max}, we have
$$\PP(A,|Z-X|\le 3, Z\in [i-1,i+1])\le ce^{-di^2}\delta,$$ and after summing over $i\in 2 \mathbb Z$ we get $\PP(A,|Z-X|\le 3)\le c\delta$.
For the second case, since we are already assuming $|X-Z|\ge 3$ as opposed to $|X-Z|\ge \sigma$, we may assume $\sigma=1$. Let $i\in 2\mathbb Z$, let $I=[i-1,i+1]$, $J=\R \setminus (i-2,i+2)$. Then we have
\begin{equation}\label{e:Ei}
\mathcal E_i:=\left\{\max_{I} \cB\in [\max_{J} \cB,\max_{J} \cB+\delta]\right\}\supseteq \{A, Z\in[i-1,i+1],|X-Z|>3\}\,.
\end{equation}
We want to bound the probability of $\mathcal E_i$ for each $i$. First, let $i=0$.
Let $\cA_2$ denote the second line of the Airy line ensemble starting from $(0,-1)$ to time $0$. Let $\tilde{\cA}_1=\cA\big|_{[-2,0]}$, $\tilde{\cA}_2=\cA\big|_{[0,2]}$. Let $\mathcal F$ be the sigma-algebra generated by $\cA', \cA_2$, $\mathcal A\big|_{ (-2,2)^c}$, and the affinely shifted bridge versions  $\tilde{\cA}_1^{[-2,0]}$ and $\tilde{\cA}_2^{[0,2]}$, see \eqref{e: affine shift}.
 
Then by the Brownian Gibbs property and Lemma \ref{l:BB}, the conditional law of $\mathcal A(0)$ given $\mathcal F$ is that of a Gaussian distribution with mean $M=(\cA(-2)+\cA(2))/2$ and variance { $2$}, conditioned to be above some $H$, where $H,M$ are $\mathcal F$-measurable.

Now, we want to bound the probability of the event $\mathcal E_0$. Note that when $\cA(-2),\cA(2)$ are fixed and the bridges joining those to $\cA(0)$ are fixed, if $\cA(0)$ changes by more than $2\delta$, the Airy process $\mathcal A$ on $[-1,1]$ changes by at least $\delta$ and so does the maximum of $\mathcal A+\mathcal A'$ on $[-1,1]$ for a fixed $\cA'$. Thus, if $M'$ denotes the maximum of $\mathcal A+\cA'$ on the set $J=(-2,2)^c$, then $M'$ is $\mathcal F$-measurable, and 
\begin{eqnarray*}
\PP(\mathcal E_0|\mathcal F)=\PP(\max_{I} \cB\in [M',M'+\delta]|\mathcal F)\le \PP(\cA(0)\in [R,R+2\delta]|\mathcal F)
\end{eqnarray*}
for some $R$ which is $\mathcal F$-measurable. We now consider the event $\cG_0=\{H-M\le \ell\}\in \cF$, with $\ell=c_3\log^{2/3}(1/\delta+1)$. Since given $\cF$, the conditional law of $\cA(0)$ is normal with variance { $2$}, mean $M$, conditioned to be at least $H$,  Lemma \ref{l:normal-delta} implies $$
\one_{\cG_0}\PP(\cE_0|\mathcal F)\le \one_{\cG_0}\PP(\cA(0)\in [R,R+2\delta]|\mathcal F)
\le \one_{\cG_0}c((H-M)\vee 1)\delta\le c\ell\delta.
$$
On the other hand, since $\cA(0) \ge H$, we have
\begin{align*}
\PP\br{\cG_0^c}&\le  \PP\left(\cA(0)-\frac{\cA(-2)+\cA(2)}{2}>\ell\right) \\&\le \PP(\cA(0)\ge \ell/2)+\PP(\cA(-2)\le -\ell/2)+\PP(\cA(2)\le -\ell/2)\le ce^{-d\ell^{3/2}}\le c \delta
\end{align*}
with the right choice of $c_3$. In conclusion, we have
$$
\PP\br{\cE_0} \le \PP\br{\cG_0^c}+\EE[\one_{\cG_0}\PP(\cE_0|\mathcal F)]\le c \delta \log^{2/3}(1+1/\delta)
.$$

For $i\neq 0$ everything above remains the same, except that we have to bound $\PP(\cA(i)-\frac{\cA(i-2)+\cA(i+2)}{2}>\ell)$  instead.  But this quantity does not depend on $i$ by the stationarity of $\cA(i)+i^2$. Thus for all $i$,
\[\PP(\mathcal E_i)\le 2C\delta\log^{2/3} (1/\delta+1)\,.\]

Finally, by Lemma $9.5$ of \cite{DOV}, for $m$ large enough, $\PP(|Z|\ge m)\le e^{-cm^3}$
for some $c>0$. Taking $m=c_4\log^{1/3}(1/\delta+1)$, for some $c_4$ large enough, by \eqref{e:Ei} we have 
\[\PP[A,|X-Z|>3]\le \sum_{i=2\mathbb Z, |i|\le m}\PP(\mathcal E_i)+\PP(|Z|>m)\le c\delta \log(1/\delta+1)\,.\qedhere\]
\end{proof}

\subsection{The main theorem on detours}\label{ss:detour-main}

We are now ready to prove the main theorem of Section \ref{sec:detour}.

\begin{theorem}[Good detours are unlikely]
    \label{t:detour}
Let $\gamma$ be a geodesic from $(x,r')$ to $(y,t')$. Let $(r,r+\sigma^3)\subset (r',t')$. 
Let $\Pi$ be the set of paths  with the same endpoints as $\gamma$ that agree with $\gamma$ on $[r',t'
]\setminus(r-2\sigma^3,r+3\sigma^3)$ and are disjoint from $\gamma$ on $(r,r+\sigma^3)$. %
Then for all { $0<\delta\le1$} we have
\begin{equation}
 \PP\br{\inf_{ \Pi}\Delta \le \delta \sigma} \le c\delta\log(1+1/\delta)^{7/3}.
\end{equation}
\end{theorem}

\begin{proof}{
By shift, scaling and shear invariance, we can assume $x=y=0$, $\sigma=1$ and $r=0$. By changing $c$, it suffices to consider sufficiently small $\delta>0$.
Let
$$\Pi_\delta=\big\{\pi\in\Pi:|\gamma|-|\pi|\le\delta\big\}.$$
Let $n=\max\{4,\lceil c_1\log(1/\delta)\rceil\}$, $t_i=i/n$ and $\chi=\eps_0n^{-2/3}$.
By Proposition \ref{p:detour-fkg}, with $\sigma'=1$, we get
\[
\PP\Big(\exists\pi\in\Pi_\delta:\max_{1\le i\le n-1}|\gamma(t_i)-\pi(t_i)|\le\chi\Big)
\le ce^{-dn}\le c\delta,
\]
if we choose $c_1\ge1/d$; here $\delta\le c_1^{\rm near}n^{2/3}$ for sufficiently small $\delta$, where $c_1^{\rm near}$ is the constant in that proposition.
Next, we apply Proposition \ref{p:detour-far} with $\sigma'=\sqrt\chi$ and $\delta'=\delta/\sigma'$. For $1\le i\le n-1$,
\[
\min\{t_i-r',t'-t_i\}\ge n^{-1},\qquad
\sigma'=\sqrt{\eps_0}\,n^{-1/3}.
\]
Thus the time-separation hypothesis holds with a universal parameter $\max\{1,\sqrt{\eps_0}\}$. After a union bound over these $n-1$ times, we get
\[
\PP\Big(\exists\pi\in\Pi_\delta,\ 1\le i\le n-1:
|\gamma(t_i)-\pi(t_i)|>\chi\Big)
\le cn\delta'\log(1+1/\delta')
\le c\delta\log(1+1/\delta)^{7/3}.
\]
This bounds $\PP(\Pi_\delta\ne\varnothing)$. Finally,
$\{\inf_{\pi\in\Pi}\Delta(\pi)\le\delta\}\subseteq\{\Pi_{2\delta}\ne\varnothing\}$.
Applying the bound with $2\delta$ and changing $c$ completes the proof.
}
\end{proof}

\section{Triple geodesics exponent}
\label{sec:triple}
The aim of this section will be to prove Theorem \ref{t: 3gd}. Let us recall the statement that we are trying to prove.
\begin{theorem}[Theorem \ref{t: 3gd}] 
    Let $A(y,x)$ be the event that there are three disjoint geodesics that start at $[-x,x]\times\{0\}$ and end in $[-y,y]\times\{1\}$. Then, for any $\eta>0$ there exists $c>0$ such that for all $x>0,y\ge 1$
    \[
        \P{A(y,x)} \leq c \br{x^{3/2 -\eta} y^3}.
    \] 
\end{theorem}
Having a bound on $A(y,x)$ will enable certain geometric constructions in the next section. The argument for Theorem \ref{t: 3gd} is similar in flavor to the argument for \cite[Proposition 4.4]{dauvergne202527}. We shall show that the event $A(y,x)$ forces the parabolic Airy line ensemble $\cA$ corresponding to $\cL(0,0; \cdot, 1)$ to behave in an uncharacteristic way due to the continuity of the extended Airy sheet, Lemma \ref{l: esc}. In particular, the top three Airy lines will be forced to be very close to each other as if they have been ``pinched''. We shall estimate the probability of this uncharacteristic behavior by comparing the Airy line ensemble with Brownian bridges using the Brownian Gibbs property (Theorem \ref{t: weiner}). 

\subsection{Non-intersecting Brownian bridges estimates}

We borrow notation from Section \ref{ss: airy}. Let $\bfx, \bfy \in \bR^k_>$ and let $p_t(\bfx, \bfy)$ denote the probability that $k$ Brownian bridges from $\bfx$ to $\bfy$ running for time $t$ never intersect. Our main tool will be the Karlin-McGregor formula for Brownian bridges (see p.48 in \cite{hammond2022brownian}), which states that
\begin{equation} \label{eq: KM}
    p_t(\bfx, \bfy) = \exp\br{-\sum_{i=1}^k x_iy_i/t}\det(e^{x_i y_j /t})_{1 \leq i,j \leq k}
\end{equation}
In particular for $k=2$, we have that
\begin{equation} \label{eq: KM2}
    p_t(\bfx, \bfy) = 1-e^{-(x_1 -x_2)(y_1-y_2)t^{-1}} \leq (x_1 -x_2)(y_1-y_2)t^{-1}.
\end{equation}
The last inequality simply uses that $1-e^{-x} \leq x$ for $x>0$. We can also make a similar estimate for three Brownian motions.
\begin{lemma} \label{l: KM3}
     For $\bfx, \bfy \in \bR^3_>$, we have that
     \[
     p_t(\bfx, \bfy) \leq \prod_{i =1,2} (x_i- x_{i+1})(y_i- y_{i+1})t^{-1} \notag.
     \]
\end{lemma}
\begin{proof}
    Note that $p_t(\bfx, \bfy)$ only depends on $t$ and the differences $x_i -x_j$ and $y_i-y_j$. Thus we can assume that $x_3 =y_3 =0$ and \eqref{eq: KM} simplifies to
    \begin{align*}
        p_t(\bfx, \bfy) &= \exp\br{-x_1y_1/t -x_2y_2/t} \begin{vmatrix}
            e^{x_1y_1/t} & e^{x_1y_2/t} & 1 \\
            e^{x_2y_1/t} & e^{x_2y_2/t} & 1 \\
            1            & 1            & 1 
        \end{vmatrix}\\
        &= \begin{vmatrix}
           1- e^{-(x_1-x_2)(y_1-y_2)/t} & 0 & e^{-x_1y_1/t}(1 - e^{(x_1-x_2)y_2/t}) \\
            e^{x_2(y_1-y_2)/t} & 1 & e^{-x_2y_2/t} \\
            1            & 1            & 1 
        \end{vmatrix}\\
        &= (1- e^{-(x_1-x_2)(y_1-y_2)/t})(1- e^{-x_2y_2/t}) %
+ e^{-x_1y_1/t}(1 - e^{(x_1-x_2)y_2/t}){ (e^{x_2(y_1-y_2)/t}-1)}.
    \end{align*}
    The second equality combines two steps of simplification. First we divide the first row by $e^{x_1y_1/t}$ and the second row by $e^{x_2y_2/t}$ to take care of the terms outside the determinant. Next we subtract the second row multiplied by $e^{-x_1(y_1-y_2)/t}$ from the first row. Notice that in the above equation, we have that
    \[
        e^{-x_1y_1/t}(1 - e^{(x_1-x_2)y_2/t}) \leq 0,
    \]
    as $x_1 >x_2 >0$ and $y_1, y_2 >0$. Thus we can conclude that
    \[
        p_t(\bfx, \bfy) \leq (1- e^{-(x_1-x_2)(y_1-y_2)/t})(1- e^{-x_2y_2/t}).
    \]
    The desired inequality now follows by applying $1-e^{-x} \leq x$ to both the terms in the right hand side.
\end{proof}

\begin{lemma} \label{l: b_bad}
Fix $ m \geq 1, \epsilon \in (0,1/100)$. Let $B$ be a $3$-tuple of Brownian bridges from $(0,\bfa)$ to $(1,\mathbf{b})$ for vectors $\bfa, \mathbf{b} \in  \bc{z \in \bR^3_> : z_1 -z_3 < m}$. For $0<\delta<1$, let $V(\delta; \bfa, \bfb)$ be the event that there exists $w \in (\epsilon, 1-\epsilon)^3_{\leq}$ such that
\begin{equation*}
    B_1, B_2, B_3 \text{ are non-intersecting},
\end{equation*}
\begin{eqnarray*}
    \abs{B_1(w_{1})- B_{2}(w_{1})} \leq \delta, \abs{B_2(w_{2})- B_{3}(w_{2})} \leq \delta, \abs{B_1(w_{3})- B_{2}(w_{3})} \leq \delta.
\end{eqnarray*}
Then for any $\eta>0$,
\[
    \P{V(\delta; \bfa, \bfb)} \leq c_{m, \eta} \delta^{3-\eta},
\]
where $c_{m, \eta}$ can be bounded by a polynomial in $m$ (for any fixed $\eta$).
\end{lemma}

\begin{proof}
    Let $ \alpha >0$ define $V^\alpha(\delta; \bfa, \bfb)$ as the event that $V(\delta; \bfa, \bfb)$ happens with $w \in  (\epsilon, 1-\epsilon)^3_{\leq} \cap \alpha \bZ^3$.
    Let $0<\eta<1/4$. The differences $B_1-B_2$ and $B_2-B_3$ are $(1/2-\eta/2)-$H\"older continuous with the H\"older constant $C_\eta$.
    Let $\beta_i(t)=B_i(t)-(1-t)a_i-tb_i$. These are Brownian bridges with zero endpoints. Let $K_\eta$ be the larger of the corresponding H\"older constants for $\beta_1-\beta_2$ and $\beta_2-\beta_3$. Since the mean gaps have slopes bounded in absolute value by $m$, we have
    \[
    C_\eta\le m+K_\eta,\qquad
    \P{K_\eta>u}\le c_\eta e^{-d_\eta u^2},
    \]
    where the law of $K_\eta$ does not depend on $\bfa,\bfb$. We may assume $\delta\le(2m)^{-2/\eta}$, since otherwise $\P{V(\delta;\bfa,\bfb)}\le1\le(2m)^{6/\eta}\delta^3$.
    Let $H=\{C_\eta\leq\delta^{-\eta/2}\}$. Then
    \[
    \P{H^c}\le\P{K_\eta>\tfrac12\delta^{-\eta/2}}
    \le c_\eta e^{-d_\eta\delta^{-\eta}},
    \]
    uniformly in $\bfa,\bfb$. This decays faster than any power of $\delta$. By a simple union bound argument it suffices to show the desired bound on this event. Notice that on this event,
    \[
        V(\delta; \bfa, \bfb) \subseteq V^{\alpha}(2 \delta; \bfa, \bfb),
    \]
    for $\alpha =\delta^{(2+\eta)/(1-\eta)}$. To bound the probability of the right hand side, we perform a union bound over all permissible values of $w$. Towards this purpose, define $V(w_0, \delta; \bfa, \bfb)$ as the event that the conditions defining $V(\delta;\bfa,\bfb)$ hold for $w=w_0$. In the bounds below, replacing $\delta$ by $2\delta$ only changes the constants.
    
    Let $g_1 =w_2-w_1$, and $g_2 = w_3-w_2$ and divide the question into two cases, when both $g_1$ and $g_2$ are non-zero and when either one of them is zero. Let us denote these events as $\mathbf{Case_1}$ and $\mathbf{Case_2}$ respectively.

    \textbf{Case 1:} Fix $w \in (\epsilon, 1-\epsilon)^3_{<} \cap \alpha \bZ^3$. Note that this forces both of them to be at least $\alpha$. Let $\bfz^i = \br{B_1(w_i), B_2(w_i), B_3(w_i)}$ for $i =1,2,3$. Thus, we have that
    \begin{align}\label{eq: c1_main}
        \P{V(w, \delta; \bfa, \bfb){\cap H}} \leq \int_{\Xi_H} p_{w_1}(\bfa, \bfz^1) & p_{g_1}(\bfz^1, \bfz^2) p_{g_2}(\bfz^2, \bfz^3) p_{1-w_3}(\bfz^3, \bfb)  d\mu(\bfz),
    \end{align}
    where $\Xi = \bc{\bfz^1, \bfz^2, \bfz^3 \in R^3_{\geq}  : z^1_1 -z^1_2 \leq \delta, z^2_2 -z^2_3 \leq \delta, z^3_1 -z^3_2 \leq \delta}$ and $\mu$ denotes the measure associated with $\bfz$. We denote by $\Xi_H$ the subset of $\Xi$ on which
    \[
    |(z_i^j-z_{i+1}^j)-(z_i^k-z_{i+1}^k)|
    \le\delta^{-\eta/2}|w_j-w_k|^{(1-\eta)/2},
    \quad i=1,2,\quad j,k=1,2,3.
    \]
    By Lemma \ref{l: KM3}, we have that the integrand is less than
    \begin{align*}
        \prod_{i \in \bc{1,2}} (a_i -a_{i+1})(z^1_i- z^1_{i+1})^2 (z^2_i-z^2_{i+1})^2 (z^3_i-z^3_{i+1})^2 (b_i-b_{i+1}) %
(w_1 g_1 g_2(1-w_3))^{-{ 2}}.
    \end{align*}
    Now since $z^2_2 -z^2_3 \leq \delta$, we can use the continuity of Brownian paths to claim that $z^1_2 -z^1_3 \leq c\delta^{-\eta/2} g_1^{1/2 - \eta/2}$ and $z^3_2 -z^3_3 \leq c\delta^{-\eta/2} g_2^{1/2 - \eta/2}$. Similarly, we can claim that $z^2_1 -z^2_2 \leq c\delta^{-\eta/2} \min(g_1,g_2)^{1/2 - \eta/2}$. We can now use these bounds to simplify the previous equation and show that the integrand is bounded by
    \begin{align*}
        c m^4 \delta^{6-3\eta} (g_1 g_2 \min(g_1, g_2))^{1-\eta}  (w_1 g_1 g_2(1-w_3))^{-2} 
        \leq c_\eps m^4 \delta^{6-3\eta} (g_1 g_2)^{-1/2-3\eta/2}.
    \end{align*}
    { This bound holds on $\Xi_H$.} Thus we can rewrite equation \eqref{eq: c1_main} as 
    \[
        \P{V(w, \delta; \bfa, \bfb){\cap H}} \leq c_\eps m^4 \delta^{6-3\eta} (g_1 g_2)^{-1/2-3\eta/2} \: \P{\bfz \in \Xi}.
    \]
    On the event $\Xi$ we can generate three independent normal random variables which are restricted to take values in a $\delta$ interval. Two of them will have a variance of at least $\epsilon$ and the last one will have a variance of order $g_1+g_2$. We prove this claim by applying the Gram-Schmidt algorithm to $\cZ_1 = z^1_1 -z^1_2, \cZ_2 = z^2_2 -z^2_3,$ and $ \cZ_3= z^3_1 -z^3_2$. The modified variables will be
    \begin{gather*}
        \tilde{\cZ_1} = \cZ_1,\;\; \tilde{\cZ_2} = \cZ_2 - \frac{\mathrm{Cov}(\tilde{\cZ_1}, \cZ_2)}{\mathrm{Var}(\tilde{\cZ_1})} \tilde{\cZ_1},\;\; 
        \tilde{\cZ_3} = \cZ_3 - \frac{\mathrm{Cov}(\tilde{\cZ_1}, \cZ_3)}{\mathrm{Var}(\tilde{\cZ_1})} \tilde{\cZ_1} -\frac{\mathrm{Cov}(\cZ_3, \tilde{\cZ_2})}{\mathrm{Var}(\tilde{\cZ_2})} \tilde{\cZ_2}.
    \end{gather*}
    These random variables will be independent and we just need to check their variances which can be evaluated as follows:
    \begin{gather*}
        \mathrm{Var}\br{\tilde{\cZ_1}} = \mathrm{Var}\br{\cZ_1}, \;\;
        \mathrm{Var}\br{\tilde{\cZ_2}} = \mathrm{Var}\br{\cZ_2} - \frac{\mathrm{Cov}(\tilde{\cZ_1}, \cZ_2)^2}{\mathrm{Var}(\tilde{\cZ_1})},\\
        \mathrm{Var}\br{\tilde{\cZ_3}} = \mathrm{Var}\br{\cZ_3} - \frac{\mathrm{Cov}(\tilde{\cZ_1}, \cZ_3)^2}{\mathrm{Var}(\tilde{\cZ_1})} -\frac{\mathrm{Cov}(\cZ_3, \tilde{\cZ_2})^2}{\mathrm{Var}(\tilde{\cZ_2})}.
    \end{gather*}
    Using the covariance formula for Brownian bridges, we have that
    \begin{align*}
        \mathrm{Var}\br{\tilde{\cZ_1}} &= 2w_1 (1-w_1) , \\
        \mathrm{Var}\br{\tilde{\cZ_2}} &= 2w_2 (1-w_2) - \frac{w_1^2(1-w_2)^2}{2w_1(1-w_1)} \geq \frac{3}{2} w_2(1-w_2),
    \end{align*}
    where the last inequality follows from $w_1 \leq w_2$. Now for the third term, we can apply the above estimate for $\tilde{\cZ_2}$,
    \begin{align*}
        \mathrm{Var}\br{\tilde{\cZ_3}} &\geq 2w_3(1-w_3) - \frac{4w_1^2 (1-w_3)^2}{2w_1(1-w_1)} - \frac{(1-w_3)^2(w_1-w_2)^2}{(1-w_1)^2 . (3/2) w_2(1-w_2)} \\
        &= \frac{2(1-w_3)(w_3-w_1)}{1-w_1} - \frac{2(1-w_3)^2(w_1-w_2)^2}{3(1-w_1)^2 w_2(1-w_2)} \\
        &\geq \frac{2(1-w_3)}{1-w_1} \bc{(w_3-w_1) - \frac{(w_1-w_2)^2}{3(1-w_1) w_2}}.
    \end{align*}
    In the last inequality we have used that $1-w_3<1-w_2$ to bound the ratio by 1. 
    \[
        \frac{(w_1-w_2)^2}{3(1-w_1) w_2} < \frac{(w_1-w_2)^2}{3(w_2-w_1)} < \frac{1}{3}(w_3-w_1),
    \]
    for $0 < w_1 < w_2<w_3< 1$. This allows us to claim that,
    \[
        \mathrm{Var}\br{\tilde{\cZ_3}} \geq \frac{4(1-w_3)(w_3-w_1)}{3(1-w_1)} \ge c_\eps (g_1+g_2).
    \]
    This proves our claim. Furthermore, since $w_1, w_2, w_3$ are bounded away from $0$ and $1$, it is a straightforward computation to see that on the event $\Xi$, $\tilde{\cZ_1}, \tilde{\cZ_2}$ and $\tilde{\cZ_3}$ will also be forced to be in an interval of length at most $2\delta$ and thus   

    \[
        \P{\Xi} \leq c_\epsilon \delta^3 (g_1+g_2)^{-1/2} \leq c_\epsilon \delta^3 g_1^{-1/4} g_2^{-1/4},
    \]
    where we use $\mathrm{AM} \geq \mathrm{GM}$ in the last inequality.
    
    All that remains is to perform the union bound over $w$ satisfying $\mathbf{Case_1}$,
    \begin{align*}
        \P{\mathbf{Case_1}{\cap H}} &\leq \sum_{w \in (\epsilon, 1-\epsilon)^3_{<} \cap \alpha \bZ^3} \P{V(w, \delta; \bfa, \bfb){\cap H}} \\
        &\leq \sum_{g_1, g_2 \in \bc{\alpha, 2 \alpha, \cdots, 1}} c_\eps m^4 \delta^{6-3\eta} (g_1 g_2)^{-1/2-3\eta/2} \delta^3 g_1^{-1/4} g_2^{-1/4} \alpha^{-1} \\
        &\leq  c_{\epsilon} m^4 \delta^{9-3\eta} \alpha^{-1 - 3 \eta} \sum_{g_1, g_2 \in \bc{\alpha, 2 \alpha, \cdots, 1}} g_1^{-3/4} g_2^{-3/4}  \\
        &\leq c_{\epsilon, m} \delta^{9-3\eta} \alpha^{-3 - 3\eta} = c_{\epsilon, m} \delta^{3-  O(\eta)}.
    \end{align*}
    In the above equation, we have claimed that the sum over $g_1, g_2$ is of the order $\alpha^{-2}$ by interpreting it as the Riemann integral of the function $f(x,y) = (xy)^{-3/4}$.

    \textbf{Case 2:} In this case there is a location $v \in (\epsilon, 1-\epsilon) \cap \alpha \bZ$ such that $B_1(v) -B_2(v) < \delta$ and $B_2(v) -B_3(v) < \delta$, where $B_i$ are still non-intersecting Brownian bridges. Let us denote this event by $U(v, \delta; \bfa, \bfb)$. We shall bound the probability of $\mathbf{Case_2}$ by applying an union bound on all possible values of $v$.

    Let $\bfz = \br{B_1(v), B_2(v), B_3(v)}$, then we can claim that in this case
    \begin{equation} \label{eq: c2_main}
        \P{U(v, \delta; \bfa, \bfb)} \leq \int_\Xi p_{v}(\bfa, \bfz) p_{1-v}(\bfz, \bfb) d\mu(\bfz),
    \end{equation}
    where  $\Xi = \bc{\bfz \in R^3_{\geq}  : z_1 -z_2 \leq \delta, z_2 -z_3 \leq \delta}$. Applying \eqref{eq: KM} to the integrand gives us that
    \[
        p_{v}(\bfa, \bfz) p_{1-v}(\bfz, \bfb) \leq c_{\epsilon} m^4 \delta^4.
    \]
    We can bound $\P{\Xi} \leq c_\eps \delta^2$ and claim that
    \[
        \P{U(v, \delta; \bfa, \bfb)} \leq c_{\epsilon} m^4 \delta^4 \: \P{\Xi} \leq c_{\epsilon} m^4 \delta^6.
    \]
    Performing the union bound for $\mathbf{Case_2}$ gives us that
    \[
        \P{\mathbf{Case_2}} \leq \sum_{v \in (\epsilon, 1-\epsilon) \cap \alpha \bZ} c_{\epsilon} m^4 \delta^6 \leq c_{\epsilon} m^4 \delta^{4-3\eta/(1-\eta)}.
    \]
    The claimed inequality now follows from combining the bounds for both cases with the bound on $\P{H^c}$, and taking $\eta$ small enough. The constants are bounded by $c_{\epsilon,\eta}(m^4+(2m)^{6/\eta})$, and hence by a polynomial in $m$.
\end{proof}

In the above lemma, we have studied one of the possible configurations that we will be encountering in the proof of Theorem \ref{t: 3gd}. There are two other configurations that will be of our interest. Namely when one drops the requirement for the third pinch $w_3$ and when one only asks for one pinch $w_1$. Similar bounds for both of these cases can be calculated using the same technique as in Lemma \ref{l: b_bad} and we state the final result as a corollary.
\begin{corollary} \label{cor: b_bad}
    In the setting of Lemma \ref{l: b_bad}, let $n\in \bc{1,2}$ and consider the event $V^n(\delta, \bfa, \bfb)$ that there exists $w \in (\eps, 1-\eps)^n_\leq$ such that
    \begin{equation*}
        B_1, B_2, B_3 \text{ are non-intersecting},
    \end{equation*}
    and for $1 \leq i \leq n$, 
    \begin{eqnarray*}
        \abs{B_i(w_{i})- B_{i+1}(w_{i})} \leq \delta.
    \end{eqnarray*}
    Then for any $\eta>0$,
    \[
        \P{V^n(\delta; \bfa, \bfb)} \leq c_{m, \eta} \delta^{n-\eta},
    \]
    where $c_{m, \eta}$ can be bounded by a polynomial in $m$ (for any fixed $\eta$).

\end{corollary}

With this, we are ready to prove Theorem \ref{t: 3gd}.
\subsection{Proof of Theorem \ref{t: 3gd}}
\begin{proof} We may assume $x\in(0,1)$ and $y>1$, or else the claim is trivial. 
    Let $(p_i, q_i)$ denote the endpoints of the three disjoint geodesics from $[-x,x]\times\{0\}$ to $[-y,y]\times\{1\}$ for $i =1,2,3$. For $y>1$, we cannot apply Corollary \ref{c: esc-unif} directly. We deploy a mesh of size 2 on $[-y,y]\times\{1\}$. Let $S_y$ denote the set of midpoints of the mesh intervals. We now apply Corollary \ref{c: esc-unif} with $z =(\mathbf{0}; \mathbf{s}), \mathbf{s} \in S_y^3$ and take a union bound to control the maximum $C$ (call it $C_y$). Thus, we can claim that
    \[
        \abs{\cL(\mathbf{0}, 0; \mathbf{q}, 1) - \cL(\mathbf{p}, 0; \mathbf{q},1)} \leq C_y x^{1/2} \log^{2/3}(c/x) +12xy, \qquad \P{C_y > a} \le cy^3e^{-ca^{3/2}}.
    \]
    By similar arguments, we can also claim that for the same $C_y$,
    \[
        \abs{\cL(0,0; q_i, 1) - \cL(p_i,0; q_i, 1)} \leq C_y x^{1/2}\log^{2/3}(c/x) +4xy,
    \]  
    for $i =1,2,3$. However, since the geodesics from $p_i$ to $q_i$ are disjoint, we have that $\cL(\mathbf{p}, 0; \mathbf{q},1) = \sum_{i=1}^3 \cL(p_i,0; q_i, 1)$. Thus, we can claim that
    \begin{equation} \label{eq: 3g(i)}
        \abs{\cL(\mathbf{0}, 0; \mathbf{q}, 1) -\sum_{i=1}^3 \cL(0,0; q_i, 1)} \leq 4 C_y x^{1/2}\log^{2/3}(c/x) +24xy.
    \end{equation}
    Consider the event $A_M := \bc{C_y \leq M \log(1/x)^{2/3}}$, then we have that
    \[
        \P{A_M} \geq 1 - cy^3e^{-cM^{3/2} \log(1/x)} \geq 1 - x^2y^3,
    \]
    for $M> M_0$. We now define $\delta = 4M_0 \log(c/x)^{4/3} x^{1/2}+ 24 xy$ and note that the L.H.S. in \eqref{eq: 3g(i)} is less than $\delta$ on the event $A_{M_0}$. 

    We now analyze the parabolic Airy line process corresponding to the extended directed landscape values, $\cL((0,0)^k, (z,1)^k)$ for $z \in [-2y, 2y]$. In particular, we will be interested in the top three lines $(\cA_1, \cA_2, \cA_3)$ and we will analyze them on the compact interval $[-y,y]$. Let us define
    \begin{multline}
        U_\delta := \big\{(f_1, f_2, f_3) \in C(\bR^3) \:\vert\: f_1 > f_2 > f_3, f_1(z_1) -f_2(z_1) \leq \delta,\\ f_2(z_2) -f_3(z_2) \leq \delta, f_1(z_3) -f_2(z_3) \leq \delta, \;\mathrm{for\; some}\; \bfz \in [-y,y]^3_< \big\}.
    \end{multline}
     The line ensemble $\cA$ is non-intersecting and equation \eqref{eq: 3g(i)} forces $(\cA_1, \cA_2, \cA_3)$ to be in $U_\delta$ following the arguments from \cite[Proof of Lemma 4.5]{dauvergne202527}. We shall abuse notation to denote this event by $\cA \in U_\delta$.

     Let us refer to $z_i$ as the pinch points. We now partition the interval $[-y-2,y+2]$ by consecutive intervals of length 1. Call these smaller intervals windows. Now, each pinch point $z_i$ must be in one of these windows, call it $W_i$. To apply our Brownian bridge estimates we need control over the Airy lines at the endpoints of the windows and we need to make sure that the pinch points are bounded away from the endpoints of their windows. 

    For the first part, we introduce a parameter $m >0$ and define the set
    \begin{equation} \label{eq: E_m}
        E_m := \bc{ (f_1, f_2, f_3) \in C(\bR^3) \:|\: \abs{f_1(i) -f_3(i)} < m, \forall\: i \in \bZ \cap [-y-2,y+2]}.
    \end{equation}
    We can apply a union bound to claim that
    \begin{align*}
        \P{\cA \in E_m^c} &\leq (2y+4) \P{\abs{\cA_1(0) -\cA_3(0)} > m} \\&\leq c y (\P{\cA_1(0) > m/2} + \P{\cA_3(0) < -m/2} )\leq c y e^{-cm^{3/2}},
    \end{align*}
    where the final inequality follows from the Tracy--Widom tail bounds of the Airy line ensemble (see {\cite[Lemma 2.1, Lemma 2.2]{dauvergne2024wiener}}).
    For the second part, we shall extend $W_i$ by appending the window to the left and to the right of it to create a window of length 3. If the window to the left also had a pinch point then we will append the next window to the left as well. Same treatment is given to the windows to the right.  Since there are only three pinch points, we will have that all pinch points are in a window of size at most 5 and all pinch points are at least a unit distance away from the boundaries. Finally, we append extra windows to make sure all extended windows containing pinch points are of the same size. Doing so might introduce some overlap between windows which contain the pinch points. If this happens, we extend our window again. Ultimately, all pinch points will be in disjoint extended windows of size at most 10. Let us denote these extended windows by $\tilde{W_i}$.

    The plan is to bound the probability of $U_\delta \cap E_m$ by taking a union bound over all possible locations of $\tilde{W_i}$. There are at most $y^3$ many ways of choosing these three windows. Let us fix one such possibility, say $\tilde{W_i} = [a_i, a_i + t]$, where $3< t< 10$ denotes the common window size and $a_i \in \bZ \cap [-y-2, y+2]$. Define the set $U_\delta(\bfa)$ as the the subset of $U_\delta$ with pinch points in $[a_i, a_i + t]$. Thus, we have that
    \[
    \P{\cA \in U_\delta \cap E_m} \leq \sum_\bfa \P{\cA \in U_\delta(\bfa) \cap E_m} \leq y^3 \max_\bfa \P{\cA \in U_\delta(\bfa) \cap E_m}.
    \]
     We now apply Theorem \ref{t: weiner} with the the set $W(a) = \cup_{i=1}^3 \llbracket1,3\rrbracket \times [a_i, a_i +t]$. This gives us a line ensemble $\cL^a$ such that
     \begin{align} \label{eq: U2L}
        \P{\cA \in U_\delta(\bfa) \cap E_m} &= \P{\cL^a \in U_\delta(\bfa) \cap E_m \:\vert\: \cL^a \in \mathrm{ORD}} \notag\\
         & \leq e^{ct^3} \P{\cL^a \in U_\delta(\bfa) \cap E_m \cap \mathrm{ORD}},
     \end{align}
    with $t \leq 10$. Since we just want an upper bound, we shall relax the restriction imposed by $E_m$ to be true only on the endpoints of the extended windows. Furthermore, we will only ask the top three lines to be ordered to keep things simple.  
    
    Let $\cF_{\mathrm{out}}$ be the $\sigma$-algebra generated by $\cL^a_i(r)$ for $(i,r)\notin W(a)$. Since the top three lines carry no ordering constraint inside $W(a)$, Theorem \ref{t: weiner} (ii), applied to each window $S_j=\llbracket 1,3\rrbracket\times[a_j,a_j+t]$ separately, shows that conditionally on $\cF_{\mathrm{out}}$ the restrictions $B^{(j)}$ of $(\cL^a_1,\cL^a_2,\cL^a_3)$ to the windows are independent, and each is a triple of independent Brownian bridges between its ($\cF_{\mathrm{out}}$-measurable, ordered) boundary values, with no conditioning.

On the event $\cL^a\in U_\delta(\bfa)\cap E_m\cap \mathrm{ORD}$, each $B^{(j)}$ is non-intersecting on its window, its endpoint spreads are less than $m$, and it contains $n_j$ of the pinch points, each at distance at least $1$ from the window boundary, where $\sum_j n_j=3$. After Brownian rescaling of the window to $[0,1]$, Lemma \ref{l: b_bad} (when $n_j=3$) and Corollary \ref{cor: b_bad} (when $n_j\in\{1,2\}$, using time reversal and the reflection $(B_1,B_2,B_3)\mapsto(-B_3,-B_2,-B_1)$ for the remaining pinch orders) bound the conditional probability of this event in window $j$ by $c_m\delta^{\,n_j-\eta/3}$, uniformly over the boundary values. By conditional independence of the window environments,
\[
    \P{\cL^a \in U_\delta(\bfa) \cap E_m \cap \mathrm{ORD} \:\vert\: \cF_{\mathrm{out}}} \leq c_m \delta^{3-\eta}.
\]

    Again, since this bound is uniform over values of $\cL^a$ outside $W(a)$, the bound also holds for the unconditioned probability. Thus, by \eqref{eq: U2L} we have that
    \[
        \P{\cA \in U_\delta(a) \cap E_m} \leq c_m \delta^{3-\eta}.
    \]
    Now, we can apply a union bound to claim that
    \begin{align*}
        \P{A(y,x)} &\leq \P{A_M^c} + \P{\cA\in E_m^c} + y^3 \max_a \P{\cA \in U_\delta(a) \cap E_m} \\
        &\leq x^2y^3 + yCe^{-cm^{3/2}} + y^3 c_m \br{M_0 \log(1/x)^{4/3} x^{1/2}+ 24 xy}^{3-\eta}.
    \end{align*}
    Notice that $\delta = O(x^{1/2} \log(1/x)^{4/3})$ for values of $x,y$ such that $x^{1/2} \log(1/x)^{4/3} > O(xy)$. When this inequality is not true, we have that $y > O(x^{-1/2} \log(1/x)^{4/3})$, which implies that
    \[
        x^{3/2-\eta}y^3 > O(x^{-\eta} \log(1/x)^{4/3}).
    \]
    Thus the claim is trivially true in this range. Hence we can assume that $\delta = O(x^{1/2} \log(1/x)^{4/3})$.
    \[
        \P{A(y,x)} \leq x^2y^3 + yCe^{-cm^{3/2}} + y^3 c_m \br{\log(1/x)^{4/3} x^{1/2}}^{3-\eta} \leq c y^3 x^{3/2-\eta}.
    \]
    The final inequality follows by optimizing over $m$ while noting that $c_m$ is a polynomial in $m$ and applying some relabelling.
\end{proof}

\section{A last passage percolation argument}
\label{sec:last-passage}

We now give a counting bound for last passage percolation in a
long-range Bernoulli environment. It will be used in the proof of
the main theorem.

\begin{proposition}\label{l:last-passage}
Fix $n\ge2$, $b\ge1$, and $p\in[0,1]$. Consider the directed graph with vertex set $\mathbb Z\times\{0,\ldots,n\}$ and all edges
\[
((j,i),(k,i+1)),\qquad j,k\in\mathbb Z,\quad 0\le i<n.
\]
Each edge has weight in $\{0,1\}$, with probability at most $p$ of having weight one. Assume that the collections of edge weights at different initial times are independent; weights within one time layer need not be independent.

Call a path $z=(z_0,\ldots,z_n)$ admissible if every integer-time partition $0=\tau_0<\cdots<\tau_\ell=n$ satisfies
\begin{equation}\label{eq:lpp-partition-constraint}
\sum_{j=1}^{\ell}\frac{|z_{\tau_j}-z_{\tau_{j-1}}|^2}{\tau_j-\tau_{j-1}}
\le b n^{1/3}\ell^{2/3}(\log(2n))^{4/3}.
\end{equation}
Let $L$ be the maximal weight of an admissible path starting in $[-n,n]\times\{0\}$ and ending at time $n$.
One may additionally restrict the terminal spatial position, since this only decreases $L$. For every integer $1\le m\le n$,
\begin{equation}\label{eq:lpp-tail}
\mathbb P(L\ge m)\le(2n+1)
\left[c b(\log(2n))^{4/3}\left(\frac nm\right)^{11/6}p\right]^m
\end{equation}
for a universal constant $c$. Set
\[
A=n b^{6/11}(\log(2n))^{8/11}p^{6/11}.
\]
There is a universal $c_0$ such that, for every $v\ge0$, with
$m_v=\lceil c_0(A+\log(2n)+v)\rceil$,
\begin{equation}\label{eq:lpp-confidence}
\mathbb P(L\ge m_v)\le n^{-3}\exp(-v).
\end{equation}
Consequently, if such environments are given for each $n$ with deterministic parameters $b_n\ge1$ and $p_n\in[0,1]$, then almost surely, for all sufficiently large $n$,
\begin{equation}\label{eq:lpp-as}
L_n\le c'\left(n b_n^{6/11}(\log(2n))^{8/11}p_n^{6/11}+\log(2n)\right).
\end{equation}
The same estimate holds for fixed $b_n=b$.
\end{proposition}
\begin{proof}
We use a counting argument. Taking the unit-time partition gives $\sum_i|z_{i+1}-z_i|^2\le bn(\log(2n))^{4/3}$. Cauchy--Schwarz therefore bounds every displacement by $\sqrt b\,n(\log(2n))^{2/3}$. Thus only finitely many paths are admissible and their maximal weight is well defined.

If $L\ge m$, some admissible path passes through $m$ edges of weight one at distinct initial times
\[
0\le t_1<\cdots<t_m<n.
\]
There are $\binom nm$ choices of these times and $2n+1$ choices of the initial position $z_0$. Write the selected edges as
$(x_j,t_j;y_j,t_j+1)$, $1\le j\le m$.
To specify their spatial positions, put
\[
g_1=t_1,\qquad g_j=t_j-t_{j-1}-1\quad(2\le j\le m),
\]
and define the signed displacements
\[
u_1=x_1-z_0,\qquad u_j=x_j-y_{j-1}\quad(2\le j\le m),
\qquad v_j=y_j-x_j.
\]
Here $g_j$ counts the unselected steps immediately preceding selected edge $j$, $u_j$ is their total displacement, and $v_j$ is the displacement along that selected edge. Once the times and $z_0$ are fixed, the vector $(u_1,v_1,\ldots,u_m,v_m)$ determines all selected endpoints recursively.

The selected endpoints determine a partition
\[
0,\ t_1,\ t_1+1,\ t_2,\ t_2+1,\ldots,\ t_m,\ t_m+1,\ n,
\]
with repeated times deleted. It has at most $2m+1$ positive-duration intervals. Apply \eqref{eq:lpp-partition-constraint} to this partition. A gap of $g_j>0$ contributes $u_j^2/g_j$, and a selected edge contributes $v_j^2$. If $g_j=0$, the endpoints coincide and $u_j=0$; that gap contributes zero. Discarding the nonnegative contribution after the last selected edge gives
\begin{equation}\label{eq:lpp-skeleton-energy}
\sum_{j=1}^m\left(\frac{u_j^2}{g_j}+v_j^2\right)\le a_m,
\qquad
 a_m=3^{2/3}b n^{1/3}m^{2/3}(\log(2n))^{4/3}.
\end{equation}
Here $2m+1\le3m$, and we refer to the sum in the first inequality as the energy of the vector. This is the gain from sampling the partition at the selected edges: the budget is of order $n^{1/3}m^{2/3}$ rather than $n$. Count every integer vector satisfying this necessary condition, even if it has no admissible completion.

Let $N(g,a_m)$ be the number of these vectors and choose $\lambda=m/a_m$. Since $a_m\ge m$, we have $0<\lambda\le1$. For each vector, the indicator that its energy is at most $a_m$ is bounded by $\exp(\lambda a_m-\lambda\,\mathrm{energy})$. Summing over the integer coordinates and factoring the resulting nonnegative sums gives
\begin{align*}
N(g,a_m)&\le \exp(\lambda a_m)
\prod_{j=1}^m
\left(\sum_u\exp(-\lambda u^2/g_j)\right)
\left(\sum_{v\in\mathbb Z}\exp(-\lambda v^2)\right).
\end{align*}
For $g_j=0$, the first sum has the single allowed term $u=0$ and equals one.
For $h>0$, monotonicity on $[0,\infty)$ gives
\[
\sum_{z\in\mathbb Z}\exp(-h z^2)
\le1+2\int_0^\infty\exp(-h x^2)\,dx
=1+\sqrt{\pi/h}.
\]
Hence, using $\lambda\le1$ and integer $g_j\ge0$, the two factors for index $j$ are bounded by
\[
c\lambda^{-1/2}(g_j+1)^{1/2},\qquad c\lambda^{-1/2},
\]
respectively. Since $\lambda a_m=m$, we conclude that
\[
N(g,a_m)\le c^m\lambda^{-m}\prod_{j=1}^m(g_j+1)^{1/2}.
\]
The time gaps telescope:
\[
\sum_{j=1}^m(g_j+1)=t_m+1\le n.
\]
The arithmetic--geometric mean inequality therefore gives
\[
\prod_{j=1}^m(g_j+1)^{1/2}\le(n/m)^{m/2},
\qquad
N(g,a_m)\le\left(\frac{ca_m}{m}\right)^m\left(\frac nm\right)^{m/2}.
\]
Substituting $a_m$, the spatial count becomes
\[
N(g,a_m)\le\left[cb(\log(2n))^{4/3}\left(\frac nm\right)^{5/6}\right]^m.
\]
This bound is uniform over all choices of selected times.

For any fixed selected times and spatial vector, the $m$ selected edges lie in different time layers. By the assumed layer independence, the probability that all their weights are one is at most $p^m$. We apply independence only to a fixed collection of edges, before taking the union over paths. A union bound now gives
\[
\mathbb P(L\ge m)\le(2n+1)\binom nm
\left(\frac{ca_m}{m}\right)^m\left(\frac nm\right)^{m/2}p^m.
\]
Using $\binom nm\le(\exp(1)n/m)^m$ and absorbing numerical constants proves \eqref{eq:lpp-tail}.

For completeness, $A^{11/6}=b n^{11/6}(\log(2n))^{4/3}p$. Choose $c_0\ge1$ sufficiently large that $c c_0^{-11/6}\le\exp(-4)$ and use $m=m_v$ from the statement. If $m_v>n$, the event $L\ge m_v$ is empty. Otherwise $m_v\ge c_0A$, and \eqref{eq:lpp-tail} implies
\[
\mathbb P(L\ge m_v)\le(2n+1)\exp(-4m_v)
\le(2n+1)(2n)^{-4}\exp(-4v)
\le n^{-3}\exp(-v).
\]
This also covers $p=0$ directly, since the bound in \eqref{eq:lpp-tail} is then zero. Thus \eqref{eq:lpp-confidence} holds uniformly in all the parameters.

Apply this estimate with $v=0$, $b=b_n$, and $p=p_n$. The probabilities are summable in $n$, so Borel--Cantelli gives $L_n<m_0$ eventually, with no independence assumption across $n$. The ceiling is absorbed in the logarithmic term, proving \eqref{eq:lpp-as}. Retaining the logarithmic term is essential when the success probability is very small.
\end{proof}

\section{Proof of the main theorem}
\label{sec:trails}

Let us bring the discussion back to finite length paths, $\pi$. Recall that we refer to such paths as {\bf roads}. Also recall the definition of $\Delta(\pi)$ from the introduction,

\[
    \Delta(\pi) = \cL(u,v) -\abs{\pi},
\] 
where $(u,v)\in \bR^4_\uparrow$ denote the endpoints of the road $\pi$.

For a pair of functions $\pi_1,\pi_2$ defined over time intervals $T_1,T_2$, we { will be interested in} the set of times at which they disagree,
\[
    T(\pi_1, \pi_2) = \{t\in T_1\cup T_2: t\notin T_1\cap T_2 \mbox{ \;\;or\;\; }\pi_1(t)\not =\pi_2(t)\}
\]
Let us define
\begin{equation}
    d_\alpha(\pi_1, \pi_2) := \inf_{\cO} \sum_{I \in \cO} \abs{I}^\alpha,
\end{equation}
where the infimum is over coverings $\cO$ of $T(\pi_1, \pi_2)$ by open intervals $I$. In an informal sense, this quantity will serve as a proxy for the $\alpha-$Hausdorff measure in our arguments.

Moreover, for a Borel measurable set $S$, let 
$$
d_{\alpha,S}(\pi_1,\pi_2)= \inf_{\cO} \sum_{I \in \cO} \abs{I}^\alpha,
$$
where the infimum is over coverings $\cO$ of $S \cap T(\pi_1, \pi_2)$ by open intervals $I$. Thus for any measurable sets $S_1,S_2$ we have 
$$d_{\alpha,S_1\cup S_2}\le d_{\alpha,S_1}+d_{\alpha,S_2}.$$ 

\begin{proposition}\label{p:da}
For $\alpha >2/3$, for all $d$ satisfying the bound
\[
{ \max\!\left\{0,\frac{1-\alpha}{\alpha-1/3}\right\}}<d<1.
\]
there exists a random $C_d>0$ so that for the $(0,0;0,1)$ geodesic $\gamma$ and every $(0,0;0,1)$-path $\pi$  we have 
$$d_\alpha(\gamma,\pi) \le C_d \Delta(\pi)^{p}, \qquad \PP(C_d>u)\le c_du^{-q},$$
where $p = d/(1+d)$ and $q= 1+d$.
\end{proposition}
\begin{proof}
It suffices to consider $0<\Delta(\pi)<\infty$: if $\Delta(\pi)=0$, then $\pi=\gamma$ by uniqueness, and if $\Delta(\pi)=\infty$, the bound is immediate. Write
\[
    T(\pi, \gamma)=\bigcup_{J\in \mathcal J_{\pi}}J,
\]
a union of disjoint open intervals. Let $\mathcal I_\pi=\{I_J:J\in\mathcal J_\pi\}$, where $I_J$ is a longest open dyadic interval contained in $J$. We write
\[
I_{i,k}=(i2^{-k},(i+1)2^{-k}),\qquad k\ge0,\quad 0\le i<2^k.
\]
If $I_J=(a,a+\ell)$, maximality gives $|J|\le4\ell$ and
$J\subset(a-2\ell,a+3\ell)$. Thus
\begin{equation}\label{e: dI}
   d_\alpha(\gamma,\pi)\le4^\alpha\:|\mathcal I_\pi|_\alpha,
\end{equation}
where $|\mathcal I_\pi|_\alpha=\sum_{I\in\mathcal I_\pi}|I|^\alpha$.

For $I=I_{i,k}$, let $\Pi_I$ be the set of $(0,0;0,1)$-paths that agree with $\gamma$ on $[0,1]\setminus((i-2)2^{-k},(i+3)2^{-k})$ and are disjoint from $\gamma$ on $I$. Let
\[
Q_I=\inf_{\pi\in\Pi_I}(|\gamma|-|\pi|)\ge0.
\]
By Theorem \ref{t:detour} and countability, almost surely $Q_I>0$ for every dyadic interval $I$.
Pick an interval $J \in \cJ_\pi$, then $\pi$ and $\gamma$ must agree on the endpoints of $J$. We can now create a new path $\pi_J$ by following $\gamma$ outside the interval $J$ and following the path $\pi$ in $J$. By our choice of $I_J$, it must be true that $\pi_J \in \Pi_{I_J}$. Thus, we can claim that
\begin{equation}    \label{e:qI}
    { \Delta(\pi)\ge \sum_{J \in \cJ_\pi}\bigl(\abs{\gamma}-\abs{\pi_J}\bigr)} \ge \sum_{I\in\mathcal I_\pi}Q _I=:Q_\pi, \qquad 
\end{equation}

Let $\eta >0$, which we shall fix later. For $\ell \in \bZ$, let $S_{\ell}$ be the sum of $\alpha$ power of the length of all open dyadic intervals $I_{i,k}$ of $[0,1]$ with $Q_I/|I|^{1/3+\eta}\in (2^{-\ell-1},2^{-\ell}]$. That is, 
$$
S_\ell=%
\sum_{k=0}^\infty \sum_{i={ 0}}^{2^k-1} \frac{1}{2^{\alpha k}} \one(Q_{I_{i,k}}2^{k(1/3+\eta)}\in (2^{-\ell-1},2^{-\ell}]).
$$
Then, by a Markov inequality argument we have that  
$$
\one(S_{\ell} \ge r2^{-d\ell} \mbox{ for some } \ell)\le 
\sum_{\ell \in \bZ} 
\frac{1}{2^{-d\ell} r} 
\sum_{k={ 0}}^\infty 
\frac{1}{2^{\alpha k}}
\sum_{i={ 0}}^{2^k-1} 
\one( Q _{I_{i,k}} \le 2^{-k(1/3+\eta)-\ell}),
$$
for any $d \in \bR.$ Taking expectations and using Theorem \ref{t:detour}, we get that
\begin{equation} \label{e:Rtail}
    \PP(R > r) \le \frac{1}{r} 
    \sum_{\ell \in \bZ} \sum_{k=0}^\infty 2^{d\ell} 2^{(1-\alpha)k} (1 \wedge c2^{-\eta k -\ell} (1+{ (\eta k +\ell)_+})^{7/3}), \qquad R=\sup_\ell S_{\ell} /2^{-d\ell}.
\end{equation}
The sum is finite as long as $0< d <1$ and $d\eta > 1-\alpha$.

Now let $S_{\pi,\ell}$ be the sum of $\alpha$ power of length of intervals $I\in \mathcal I_\pi$ with $Q_I/|I|^{1/3+\eta} \in (2^{-\ell-1},2^{-\ell}]$. Then
$$
|\mathcal I_\pi|_\alpha =\sum_{\ell \in \bZ} S_{\pi,\ell},
$$
and we can claim that
\begin{align*}
    S_{\pi, \ell} &= \sum_{I \in \mathcal{I}_\pi} \one \br{\frac{Q_I}{\abs{I}^{1/3 +\eta}} \in (2^{-\ell-1}, 2^{-\ell}]} \abs{I}^\alpha \leq \sum_{I \in \mathcal{I}_\pi} \frac{Q_I}{\abs{I}^{1/3 +\eta} 2^{-\ell-1}} \abs{I}^\alpha  \leq Q_\pi 2^{\ell+1} ,
\end{align*}
where the final inequality is true for $\alpha \ge 1/3 +\eta$ as $\abs{I} \leq 1$. Combining this inequality with the restriction from \eqref{e:Rtail}, the optimal choice for $\eta$ is $\alpha-1/3$ which gives us the restriction that $\alpha > 2/3$. Thus
$$
|\mathcal I_\pi|_\alpha = \sum_{\ell\le \ell_0} S_{\pi,\ell}
+\sum_{\ell > \ell_0} S_{\pi,\ell}\le { 2^{2+\ell_0} Q_\pi} +\sum_{\ell > \ell_0} S_{\ell}\le { 2^{2+\ell_0} Q_\pi} + cR2^{-d\ell_0}.
$$
We choose $\ell_0\in\mathbb Z$ so that
\[
2^{\ell_0}\le\br{\frac{R}{Q_\pi}}^{1/(1+d)}<2^{\ell_0+1}.
\]
Using this choice of $\ell_0$ along with \eqref{e:qI}, \eqref{e: dI} we have that  
\[
d_\alpha(\gamma, \pi) \le cR^{1/(1+d)}\Delta(\pi)^{d/(1+d)}.
\]
The tail bound \eqref{e:Rtail} on $R$ implies the claim. 
\end{proof}

\begin{proposition} \label{p :cpt_do}
    For every $\alpha, d$ satisfying the bounds in Proposition \ref{p:da}, $0< \delta < 1/2,$  there exists a random $C >0$ such that for all finite length paths $\pi$ with endpoints $u \in B = [-1,1] \times \bc{0}$ and $v \in T = [-1,1] \times \bc{1}$, there exists a geodesic path $\gamma$ from $B$ to $T$ such that
    \[
        d_{\alpha, (\delta, 1-\delta)} (\gamma, \pi) \leq C \Delta(\pi)^{p}, \qquad \P{C >u} \leq c(1+\delta^{-{ 3}})u^{-q'},
    \]
    where the inequality holds for $p =d/(1+d)$ and all $q' <(1+d)/5$.
\end{proposition}
\begin{proof}
    We start by limiting the locations that geodesics across $B$ and $T$ can visit. To do so, we simply ask that the geodesics from $(\pm  1,0)$ to $(\pm 1,1)$ have transversal fluctuations less than $M \delta^{2/3}$ at times $\delta$ and $1- \delta$. { Let us call} this event $\cT$, then we have that
    \[
        \P{\cT} \geq 1 - \exp(- c M^3).  
    \]
    Now choose $\epsilon >0$ and divide $B, T$ into an $\epsilon-$mesh, call the collection of $\epsilon-$ intervals as $m(B), m(T)$. Let $\cA^B_\epsilon$ be the event that there are three disjoint geodesic paths from an interval $I \in m(B)$, going to the set $[-M\delta^{2/3}-1, M\delta^{2/3}+1] \times \delta$. We define $\cA^T_\epsilon$ as the analogue of $\cA^B_\epsilon$ for the top side and define $\cA_\epsilon = \cA^B_\epsilon \cup \cA^T_\epsilon$. We can use the scaling properties of the directed landscape to apply Theorem \ref{t: 3gd}, and perform { a union bound} to claim that for $M \delta^{2/3}< \eps^{-1/2}$,
    \[
        \P{\cA_\epsilon} \leq  c \epsilon^{1/2 - \eta} M^3 \delta^{-3} ,
    \]
    for any $\eta >0$. We shall see that at the end, the restriction on $M$ will not matter.
    
    Let $I^*, J^*$ denote the intervals in $m(B), m(T)$ which contain the endpoints $u$ and $v$. For ease of notation, we relabel the endpoints of $I^*$ as $\ell$ and $r$, denoting the left and the right endpoint. Similarly we define $\ell'$ and $r'$ for $J^*$. We modify the path $\pi$ by replacing the initial $\delta$ segment by the geodesic from $u$ to $(\pi(\delta), \delta)$. We do the same modification for the final $\delta$ segment of $\pi$ and call the modified path $\tilde{\pi}$. Note that, since we replaced parts of $\pi$ by geodesic segments, we must have that $\abs{\pi} \leq \abs{\tilde \pi}.$

    We can now claim that $\tilde{\pi}$ must intersect either $\gamma_{\ell,\ell'}$ or $\gamma_{r,r'}$ at least once before time $\delta$ and after time $1- \delta$ on the event $\cA_\epsilon^c \cap \cT$. Two cases can arise now, firstly $\tilde{\pi}$ can intersect the same geodesic in both the time intervals or they can be different. 
    \begin{figure}[tb]
        \centering
        \includegraphics[width=0.5\linewidth]{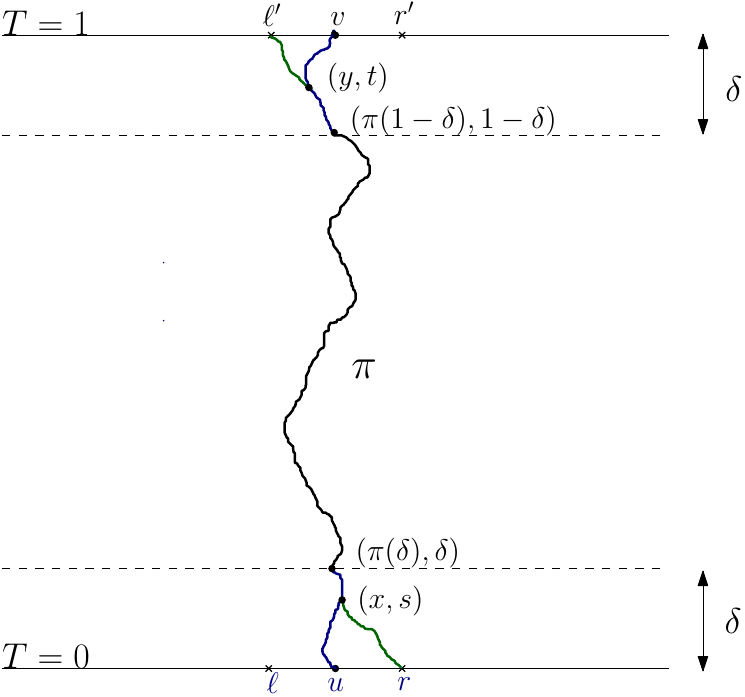}
        \caption{The black path represents the path $\pi$. We change the initial and final segments of $\pi$ by geodesic paths depicted in blue to get $\tilde \pi$ . The green paths depict segments of $\gamma_{\ell, \ell'}$ and we construct $\pi^*$ by following the green paths till we hit the blue path and then using $\tilde{\pi}$ in between the first and last intersection points.}
        \label{fig: cpt_intvl}
    \end{figure}
    
    We deal with the first case, and assume that $\tilde{\pi}$ intersects $\gamma_{\ell,\ell'}$ in both time intervals. Let us denote the first intersection point by $(x,s)$ and the final intersection point by $(y,t)$. We construct a new path $\pi^*$ starting at $\ell$, following $\gamma_{\ell, \ell'}$ till time $s$ and then following $\tilde{\pi}$ till time $t$, from where we follow $\gamma_{\ell, \ell'}$ to reach $\ell'$. The construction is depicted in Figure \ref{fig: cpt_intvl}. We can now apply Proposition \ref{p:da}, to $\gamma_{\ell,\ell'}$ and $\pi^*$ to claim that 
    \begin{align*}
        d_{\alpha, (\delta, 1-\delta)}(\gamma_{\ell,\ell'}, \pi^*) &\leq C_{\ell,\ell'}(\cL_{\ell,\ell'} - \abs{\pi^*})^{p} \\
        &=  C_{\ell,\ell'}(\cL_{(x,s),(y,t)} - \abs{\pi^*}_{[s,t]})^{p} \\
        &=  C_{\ell,\ell'}(\cL_{u,(x,s)} + \cL_{(x,s),(y,t)} + \cL_{(y,t),v} - \abs{\tilde{\pi}}  )^{p} \\
        &\leq C_{\ell,\ell'}(\cL_{u,v} - \abs{\pi})^{p}.
    \end{align*}
    The second last equality follows from noting that the initial and final segments of $\tilde{\pi}$ are geodesics from $u$ to $(x,s)$ and $(y,t)$ to $v$. Since, $\pi$ and $\pi^*$ differ only on $[0,1] \setminus [\delta,1-\delta]$, we have that
    \[
        d_{\alpha, (\delta, 1-\delta)}(\gamma_{\ell,\ell'}, \pi) =  d_{\alpha, (\delta, 1-\delta)}(\gamma_{\ell,\ell'} , \pi^*) \leq C_{\ell,\ell'}\Delta(\pi)^{p}.
    \]
    
    Now, for the other case assume that the initial intersection happens with $\gamma_{\ell,\ell'}$ while the final intersection is with $\gamma_{r,r'}$. Our claim is that $\tilde{\pi}$ will intersect $\gamma_{\ell,r'}$ at least once before time $\delta$ and after time $1- \delta$. This simply follows from the ordering of geodesic paths in these time intervals. Notice that $\gamma_{\ell,r'}$ is to the right of $\gamma_{\ell,\ell'}$, thus it must intersect $\tilde \pi$ in $(0, \delta)$. Similarly, for the top side. We can now run the same arguments with $\gamma_{\ell,r'}$ instead of $\gamma_{\ell,\ell'}$ to claim that
    \[
        d_{\alpha, (\delta, 1-\delta)}(\gamma_{\ell,r'}, \pi) \leq C_{\ell,r'}\Delta(\pi)^{p}.
    \]
    Now define $C_\eps:= \max C_{u_\epsilon, v_\epsilon}$, where the maxima is taken over all mesh points $u_\epsilon, v_\epsilon$. We can now claim that on the event $\cA_\epsilon^c \cap \cT$, there exists a geodesic from $B$ to $T$ such that 
    \[
        d_{\alpha, (\delta, 1-\delta)}(\gamma, \pi) \leq C_\eps \Delta(\pi)^{p}.  
    \]
{ Let us now take $\epsilon_j=2^{-j}$ and set $C_j=\max\{1,C_{\epsilon_j}\}$. Let
\[
G_{j,M}=\mathcal T_M\cap\mathcal A_{\epsilon_j,M}^c,
\qquad M\in\mathbb N.
\]
On $G_{j,M}$, the construction above applies to every path from $B$ to $T$. We now define
\[
C=2\inf\{C_j:G_{j,M}\text{ holds and }
M\delta^{2/3}<\epsilon_j^{-1/2}\}.
\]
All the events and constants in this countable infimum are measurable with respect to the landscape in the strip. Thus $C$ is also measurable with respect to this landscape. Taking $M=j$, the condition $M\delta^{2/3}<\epsilon_j^{-1/2}$ holds for all sufficiently large $j$, and we have
\[
\sum_j\mathbb P(G_{j,j}^c)
\le\sum_j\left(e^{-cj^3}
+c\delta^{-3}j^3 2^{-j(1/2-\eta)}\right)<\infty.
\]
By Borel--Cantelli, $G_{j,j}$ holds for all sufficiently large $j$, so $C<\infty$ almost surely. We can therefore choose a pair $(j,M)$ in the infimum with $C_j<C$. Applying the construction above with this mesh proves the required inequality for every path. If $\Delta(\pi)=0$, we take $\gamma=\pi$.

It remains to bound the tail of $C$. Fix $j,M$ satisfying $M\delta^{2/3}<\epsilon_j^{-1/2}$ and let $u\ge1$. On $G_{j,M}\cap\{C_j\le u\}$, we have $C\le2u$. Thus
\begin{align*}
\mathbb P(C>2u)
&\le\mathbb P(G_{j,M}^c)+\mathbb P(C_j>u)\\
&\le e^{-cM^3}
+c\delta^{-3}M^3\epsilon_j^{1/2-\eta}
+c\epsilon_j^{-2}u^{-q},\qquad q=1+d.
\end{align*}
For $q'<q/5$, we can choose $\eta>0$ small enough and $b$ such that
\[
\frac{q'}{1/2-\eta}<b<\frac{q-q'}2.
\]
We now choose $\epsilon_j\in(\tfrac12u^{-b},u^{-b}]$ and
$M=\lceil c_1(\log u)^{1/3}\rceil$. For all sufficiently large $u$, the condition $M\delta^{2/3}<\epsilon_j^{-1/2}$ holds for every $\delta<1/2$. Thus, choosing $c_1$ large enough, the previous bound gives
\[
\mathbb P(C>2u)\le c(1+\delta^{-3})u^{-q'}.
\]
Changing the constant $c$ gives the claimed tail bound for all $u>0$.}
\end{proof}

{
\subsection{Proof of the main theorem}
\label{ss:main-proof}

We take $k\in\mathbb N$, $n=2^{3k}$, $r=n^{-1}$, and $0<\delta\le1/4$. We write
\[
T_i=[ir,(i+1)r],\qquad
T_i^\delta=ir+[\delta r,(1-\delta)r],\qquad
T^\delta=\bigcup_{i=0}^{n-1}T_i^\delta.
\]
We fix $\alpha\in(2/3,1)$ and admissible $p,q'$ from Proposition \ref{p :cpt_do}, and write
\[
s=\frac{6q'}{11},\qquad a=ps,\qquad
\rho=\alpha-\frac{1+a/3}{1+a}.
\]
Throughout this subsection, $\alpha,p,q',\delta$ are fixed. We shall write $K_0(b,h)$ for a finite random integer which is the same for all paths in $\mathcal T_{b,h}$, and increase it when needed.

Let us consider the directed graph with vertex set
\[
V_k=\{I_{j,i}:j\in\mathbb Z,\ 0\le i\le n\},\qquad
I_{j,i}=r^{2/3}[j,j+1)\times\{ir\},
\]
and edge set
\[
E_k=\{(I_{j,i},I_{\ell,i+1}):j,\ell\in\mathbb Z,\ 0\le i<n\}.
\]
For an integer path $z=(z_0,\ldots,z_n)$, let
\[
e_i(z)=(I_{z_i,i},I_{z_{i+1},i+1}),\qquad 0\le i<n.
\]
For an edge $e=(I,J)$ between times $ir$ and $(i+1)r$, let $I^+,J^+$ be the concentric closed intervals of length $2r^{2/3}$. Let $C_e\ge1$ be the constant from Proposition \ref{p :cpt_do} applied to $I^+,J^+$, after KPZ scaling and shear invariance. This constant depends only on the landscape between these times, and we have
\begin{align}
\mathbb P(C_e>u)&\le c_\delta u^{-q'},\qquad u\ge1,
\label{eq:local-constant-tail}\\
\forall\pi':I\to J,\ \Delta(\pi')>0,\quad
\exists\gamma\in\Gamma(I^+,J^+):\quad
 d_{\alpha,T_i^\delta}(\gamma,\pi')
&\le C_e r^{\alpha-p/3}\Delta(\pi')^p.
\label{eq:strip-comparison}
\end{align}
The collections $(C_e:e\text{ in layer }i)$ are independent for different $i$.

For each deterministic integer $m\ge1$, let $\mathcal P_{k,m}$ be the set of integer paths $z=(z_0,\ldots,z_n)$ satisfying
\[
|z_0|\le n,\qquad \max_{0\le i\le n}|z_i|\le mn^{2/3}+1,
\]
and, for every integer-time partition $0=\tau_0<\cdots<\tau_\ell=n$,
\begin{equation}\label{eq:one-cutoff-path-class}
\sum_{j=1}^{\ell}
\frac{|z_{\tau_j}-z_{\tau_{j-1}}|^2}{\tau_j-\tau_{j-1}}
\le mn^{1/3}\ell^{2/3}(\log(2n))^{4/3}.
\end{equation}
Notice that $\mathcal P_{k,m}$ is deterministic, and that we use the same $m$ in both bounds.

We also use the translated cells
\[
\widehat T_i=T_i+r/2,\qquad
\widehat T_i^\delta=T_i^\delta+r/2,\qquad
\widehat T^\delta=\bigcup_{i=0}^{n-2}\widehat T_i^\delta.
\]
The translated graph has layers at $(i+1/2)r$, $0\le i\le n-1$, and the same spatial mesh. We define $\widehat e_i$ and $\widehat C_e$ as above, and define $\widehat{\mathcal P}_{k,m}$ by replacing $n$ with $n-1$ in the definition of $\mathcal P_{k,m}$.

}
\begin{lemma}\label{l:trail-grid-admissibility}{

Almost surely, for every $b,h\ge1$, set
\[
m_{b,h}=\left\lceil c(1+b+K_b+h)\right\rceil,
\qquad
\mathcal T_{b,h}=\left\{\pi:\sup_{t\in[0,1]}|\pi(t)|\le b,\quad |\pi|\ge-h\right\},
\]
where $K_b$ is the random constant from Proposition \ref{p: Mink}, depending only on the landscape and $b$, and $c$ is a universal constant. For all $k\ge K_0(b,h)$, uniformly over $\pi\in\mathcal T_{b,h}$, the path
\[
z_i=\lfloor n^{2/3}\pi(i/n)\rfloor,\qquad 0\le i\le n,
\]
belongs to $\mathcal P_{k,m_{b,h}}$. The translated path
$\widehat z_i=\lfloor n^{2/3}\pi((i+1/2)/n)\rfloor$, $0\le i\le n-1$, belongs to $\widehat{\mathcal P}_{k,m_{b,h}}$ for the same $K_0(b,h)$.
}
\end{lemma}
\begin{proof}{

Let $0=\tau_0<\cdots<\tau_\ell=n$ be a partition and write $g_j=\tau_j-\tau_{j-1}$. Applying Proposition \ref{p: Mink} at the times $\tau_j/n$, we get
\[
\sum_j\frac{|n^{2/3}\pi(\tau_j/n)-n^{2/3}\pi(\tau_{j-1}/n)|^2}{g_j}
\le hn^{1/3}+K_b\sum_jg_j^{1/3}\bigl(\log(en/g_j)\bigr)^{4/3}.
\]
Since $g_j\ge1$ and $\sum_jg_j=n$,
\[
\sum_jg_j^{1/3}\le n^{1/3}\ell^{2/3},\qquad
\sum_jg_j^{-1}\le\ell\le n^{1/3}\ell^{2/3}.
\]
By the definition of $z_i$, each displacement changes by at most one when we take integer parts. Thus
\begin{align*}
\sum_j\frac{|z_{\tau_j}-z_{\tau_{j-1}}|^2}{g_j}
&\le2hn^{1/3}+cK_b n^{1/3}\ell^{2/3}(\log(2n))^{4/3}+2\ell\\
&\le c(1+K_b+h)n^{1/3}\ell^{2/3}(\log(2n))^{4/3}.
\end{align*}
For the translated grid, we adjoin $0$ and $1$ to a partition of the shifted mesh times and apply the same bound. The two extra terms in $\sum_jg_j^{1/3}(\log(en/g_j))^{4/3}$ are at most $c(\log(2n))^{4/3}$. Since $n-1\ge n/2$, we obtain the bound defining $\widehat{\mathcal P}_{k,m_{b,h}}$ after increasing $c$.

The spatial bound gives $|z_i|,|\widehat z_i|\le bn^{2/3}+1$. We increase $K_0(b,h)$ so that $bn^{2/3}+1\le n-1$ for all $k\ge K_0(b,h)$. Choosing $c$ large enough in the definition of $m_{b,h}$ gives both spatial bounds, and proves the claim uniformly over $\mathcal T_{b,h}$.
}
\end{proof}

\begin{lemma}\label{l:uniform-threshold-count}{

Almost surely, simultaneously for every integer $m\ge1$, there is a finite random integer $K_0(m)$ such that, for all $k\ge K_0(m)$ and all $u>0$,
\begin{equation}\label{eq:uniform-threshold-LP}
L_{k,m}(u):=\max_{z\in\mathcal P_{k,m}}
\sum_{i=0}^{n-1}\mathbf1\{C_{e_i(z)}>u\}
\le c_{m,\delta}(knu^{-s}+k).
\end{equation}
The same estimate holds for
$\widehat L_{k,m}(u)=\max_{z\in\widehat{\mathcal P}_{k,m}}\sum_{i=0}^{n-2}\mathbf1\{\widehat C_{\widehat e_i(z)}>u\}$
with the same $K_0(m)$. Here $c_{m,\delta}$ is deterministic and may also depend on the fixed parameters $\alpha,p,q'$.
}
\end{lemma}
\begin{proof}{

We fix an integer $m$. For $u=2^j$, we apply Proposition \ref{l:last-passage} with
\[
p_j=1\wedge c_\delta2^{-q'j},\qquad b_n=m.
\]
The paths in $\mathcal P_{k,m}$ satisfy the condition in that proposition by \eqref{eq:one-cutoff-path-class}. The restriction on their positions can only decrease the maximum. Since $s=6q'/11$ and $(\log(2n))^{8/11}\le ck$, we get
\[
\mathbb P\!\left(L_{k,m}(2^j)>c_{m,\delta}[kn2^{-sj}+k]\right)\le n^{-3}.
\]
We choose a constant $A$ such that $Aq'>7$, and by a union bound we have that
\begin{equation} \label{eq: fail1}
    \mathbb P\!\left(\exists j\in\{0,\ldots,\lceil Ak\rceil\}:\, L_{k,m}(2^j)>c_{m,\delta}[kn2^{-sj}+k]\right) \le(\lceil Ak\rceil+1)2^{-9k}.
\end{equation}
Let $E_{k,m}$ be the set of edges whose two endpoint indices have absolute value at most $mn^{2/3}+1$. These are the edges that can possibly be used by a path in $\cP_{k,m}$. The number of such edges is at most
\[
|E_{k,m}|\le n(2mn^{2/3}+3)^2\le c_m2^{7k}.
\]
Consequently,
\begin{align} \label{eq: fail2}
\mathbb P\!\left(\max_{e\in E_{k,m}}C_e>2^{Ak}\right)
&\le c_{m,\delta}2^{(7-Aq')k}.
\end{align}
For each $k$, the same argument applies to the translated environment with $n-1$ layers. Since $n-1\ge n/2$, its two failure probabilities (\eqref{eq: fail1}, \eqref{eq: fail2}) satisfy the same bounds up to a constant. Summing the four bounds over $k$ and applying Borel--Cantelli gives a finite random $K_0(m)$ such that none of these events occurs for $k\ge K_0(m)$. For all $k\ge K_0(m)$, the monotonicity of $L_{k,m}(u)$ in $u$ gives \eqref{eq:uniform-threshold-LP} for $1\le u\le2^{Ak}$, after enlarging the deterministic constant. For $u>2^{Ak}$, we have $L_{k,m}(u)=0$, while for $u<1$ we can simply use $L_{k,m}(u)\le n$. The same monotonicity argument applies to $\widehat L_{k,m}$. Taking a countable intersection, we get both bounds for all integers $m$ on the same probability-one event.
}
\end{proof}

\begin{proposition}\label{p:subdivision}{

If $\rho>0$, then almost surely, for every $b,h\ge1$ and $0<\zeta<3\rho$, there is a finite random $K_0(b,h)$ such that, for every $\pi\in\mathcal T_{b,h}$ and every $k\ge K_0(b,h)$, there are geodesics $\gamma_{k,i}$ across $T_i$ satisfying
\[
d_{\alpha,T^\delta}(\pi,\gamma_k)
\le(1+\Delta(\pi))2^{-\zeta k},\qquad
\gamma_k=\bigsqcup_i\gamma_{k,i}.
\]
The same conclusion holds for $\widehat T^\delta$ and geodesics across the translated cells. The same $K_0(b,h)$ works for both grids and every path in $\mathcal T_{b,h}$. Notice that every road on $[0,1]$ belongs to such a class for some integers $b,h$.
}
\end{proposition}
\begin{proof}{

We fix $b,h\ge1$ and set $m=m_{b,h}$. By Lemma \ref{l:trail-grid-admissibility}, the integer path associated with any $\pi\in\mathcal T_{b,h}$ belongs to $\mathcal P_{k,m}$ for $k\ge K_0(b,h)$. Lemma \ref{l:uniform-threshold-count} holds for all integers $m$ on the same probability-one event, so we can apply it to this random choice of $m$. Increasing $K_0(b,h)$ to be at least $K_0(m_{b,h})$, both lemmas apply for every $k\ge K_0(b,h)$.

We now fix $0<\zeta<3\rho$, $k\ge K_0(b,h)$ and a road $\pi\in\mathcal T_{b,h}$. Let
\[
D=\Delta(\pi),\qquad D_i=\Delta(\pi|_{T_i}),\qquad
\sum_iD_i\le D.
\]
For each $i$, we choose $\gamma_{k,i}$ as in \eqref{eq:strip-comparison}. When $D_i=0$, we take $\gamma_{k,i}=\pi|_{T_i}$. Let
\[
h_i=d_{\alpha,T_i^\delta}(\pi,\gamma_{k,i}).
\]
Covering the disagreement set by one interval, and using \eqref{eq:strip-comparison}, we get
\begin{equation}\label{eq:local-covering-cost}
0\le h_i\le r^\alpha,\qquad
h_i\le C_{e_i(z)}r^{\alpha-p/3}D_i^p.
\end{equation}
If $D=0$, then all $h_i=0$ and there is nothing to prove. We assume $D>0$ and let $N_k(t)=\#\{i:h_i>t\}$. For any $t,u>0$, we have
\[
\{i:h_i>t\}\subseteq
\{i:C_{e_i(z)}>u\}\ \cup\
\left\{i:D_i>\left(\frac{t}{u r^{\alpha-p/3}}\right)^{1/p}\right\}.
\]
We can now apply \eqref{eq:uniform-threshold-LP} and use $\sum_iD_i\le D$ to claim that
\begin{align*}
N_k(t)&\le c_{m,\delta}(knu^{-s}+k)
+D r^{\alpha/p-1/3}t^{-1/p}u^{1/p}\\
&=x(t)u^{1/p}+y_ku^{-s}+c_{m,\delta}k,
\end{align*}
where
\[
x(t)=D r^{\alpha/p-1/3}t^{-1/p},\qquad y_k=c_{m,\delta}kn.
\]
Since $a=ps$, the two terms are equal when we take
\[
u=\left(\frac{y_k}{x(t)}\right)^{p/(1+a)}.
\]
Then
\[
x(t)u^{1/p}=y_ku^{-s}=x(t)^{a/(1+a)}y_k^{1/(1+a)},
\]
and therefore
\[
N_k(t)\le2D^{a/(1+a)}
r^{(\alpha/p-1/3)a/(1+a)}y_k^{1/(1+a)}t^{-\nu}
+c_{m,\delta}k,\qquad \nu=\frac{s}{1+a}<1.
\]
By \eqref{eq:local-covering-cost},
\[
\sum_i h_i=\int_0^{r^\alpha}N_k(t)\,dt,\qquad
\int_0^{r^\alpha}t^{-\nu}\,dt
=\frac{r^{\alpha(1-\nu)}}{1-\nu}.
\]
Using $n=r^{-1}$, the exponent of $r$ in the resulting bound is
\[
\frac{(\alpha/p-1/3)a-1}{1+a}
+\alpha\left(1-\frac{s}{1+a}\right)
=\alpha-\frac{1+a/3}{1+a}=\rho.
\]
Thus
\begin{equation}\label{eq:quadratic-content-cost}
\sum_i h_i\le c_{m,\delta}D^{a/(1+a)}k^{1/(1+a)}r^\rho
+c_{m,\delta}kr^\alpha.
\end{equation}
Finally,
\[
d_{\alpha,T^\delta}(\pi,\gamma_k)\le\sum_i h_i,\qquad
r^\rho=2^{-3\rho k},\qquad \rho<\alpha.
\]
Since $D^{a/(1+a)}\le1+D$, we get
\[
d_{\alpha,T^\delta}(\pi,\gamma_k)
\le c_{m,\delta}(1+D)
\left(k^{1/(1+a)}2^{-3\rho k}+k2^{-3\alpha k}\right).
\]
Since $\zeta<3\rho<3\alpha$, we increase $K_0(b,h)$ so that the claim holds for $k\ge K_0(b,h)$.

For the translated cells, we have $\sum_i\Delta(\pi|_{\widehat T_i})\le D$, by adjoining the two boundary intervals. The two lemmas hold for this grid with the same cutoffs, so the calculation above applies unchanged. Increasing $K_0(b,h)$ gives both bounds, uniformly over $\mathcal T_{b,h}$.
}
\end{proof}

\begin{corollary}\label{cor: apx_fnl}
{ Under the same assumptions, put
\[
E_\pi=\{t\in[0,1]:(\pi(t),t)\notin\mathcal G\}.
\]
For every $0<\zeta<3\rho$, almost surely, for each $b,h\ge1$ there is a finite random $K_0(b,h)$ such that, for every $\pi\in\mathcal T_{b,h}$ and every $k\ge K_0(b,h)$, 
\[
{ \mathcal H^\alpha_{2r}}(E_\pi)
\le(1+\Delta(\pi))2^{-\zeta k}.
\]}
\end{corollary}
\vspace{-\baselineskip}
\begin{proof}
{ We take $\delta=1/4$ and apply Proposition \ref{p:subdivision} to the grid $r\mathbb Z$ and its translate by $r/2$. The middle halves of their cells contained in $[0,1]$ cover $[r,1-r]$. On each middle half, $E_\pi$ is contained in the disagreement set with the corresponding geodesic. Restricting the covers to intervals of diameter less than $2r$, and adding the two boundary intervals at cost at most $cr^\alpha$, we get, for $\zeta<\zeta_1<3\rho$,
\[
\mathcal H^\alpha_{2r}(E_\pi)
\le2(1+\Delta(\pi))2^{-\zeta_1k}+cr^\alpha
\le(1+\Delta(\pi))2^{-\zeta k},\qquad k\ge K_0(b,h).
\]
The proposition gives a common cutoff for both grids; we increase it to obtain the last inequality.}
\end{proof}

\begin{proof}[Proof of Theorem \ref{t: main}]
{ By the corollary, we have
\[
\rho>0\quad\Longrightarrow\quad
\mathcal H^\alpha(E_\pi)=0,
\qquad
\rho>0\ \Longleftrightarrow\
\alpha>\frac{3+a}{3(1+a)}.
\]
We can choose a countable sequence of parameters satisfying the required bounds such that
\[
d\uparrow1,\qquad
p=\frac d{1+d},\qquad q'\uparrow\frac{1+d}5.
\]
Then
\[
pq'\longrightarrow\frac15,\qquad
a=\frac{6pq'}{11}\longrightarrow\frac6{55},\qquad
\frac{3+a}{3(1+a)}\longrightarrow\frac{57}{61}.
\]
For each $\alpha_j$ in a countable sequence decreasing to $57/61$, we choose these parameters so that $\rho_j>0$. On} the { intersection} of { the} { resulting probability-one events,
\[
\forall\pi:[0,1]\to\mathbb R\text{ a road},\quad
\mathcal H^{\alpha_j}(E_\pi)=0\ \ (j\ge1)
\quad\Longrightarrow\quad
\dim_H E_\pi\le\frac{57}{61}.
\]
For a general time interval, we apply the scaled statement on all rational subintervals. The claim then follows by countable stability of Hausdorff dimension, since the two endpoints do not affect the bound.}
\end{proof}

\subsection{Roads don't meet fixed paths often: the proof of Corollary \ref{c: f_int}} 

\label{s: f_int}

We will need the following folklore lemma.

\begin{lemma}\label{lem:zero-set-dimension}
Let $I$ be a compact interval and $0<\alpha\le 1$.
Suppose that $F:I\to\mathbb R$ is a random  $\alpha$-H\"older continuous function, and that the real random variable
$F(t)$ has a density $p_t$ satisfying
\[
    s=\sup_{t\in I}\|p_t\|_\infty<\infty.
\]
Then, almost surely,
\[
    \overline{\dim}_{M} F^{-1}(0)\le 1-\alpha.
\]
\end{lemma}

\begin{proof}
By rescaling, we may assume that $I=[0,1]$. Let $c$ be the
uniform density bound, and let $K$ be the $\alpha$-H\"older
seminorm of $F$, which is finite almost surely. For each
positive integer $m$, define the random closed set
\[
    Z_m=
    \begin{cases}
        F^{-1}(0), & K\le m,\\
        \varnothing, & K>m.
    \end{cases}
\]
Let $N_n(Z_m)$ count the dyadic intervals of length $2^{-n}$
meeting $Z_m$. If such an interval meets $Z_m$, its left
endpoint $j2^{-n}$ satisfies
$|F(j2^{-n})|\le m2^{-n\alpha}$. Consequently,
\[
    \mathbb E N_n(Z_m)
    \le \sum_{j=0}^{2^n-1}
        \mathbb P\bigl(|F(j2^{-n})|\le m2^{-n\alpha}\bigr)
    \le 2sm\,2^{n(1-\alpha)}.
\]
By \cite[Lemma~6.10.5]{BishopPeres2017}, an estimate
$\mathbb E N_n(A)\le c'2^{n\beta}$ for a random closed set
$A\subset[0,1]$ implies
$\overline{\dim}_{M}A\le\beta$ almost surely.
Applying this for a sequence $\beta\downarrow1-\alpha$
gives $\overline{\dim}_{M}Z_m\le1-\alpha$ almost surely.
Almost surely these bounds hold for every $m$, and $K<\infty$,
so $F^{-1}(0)=Z_m$ for some $m$. This proves the claim.
\end{proof}

The next lemma is based on results of \cite{ganguly2022fractal}.

\begin{lemma}\label{l:GZ}
Let $\gamma=\gamma_{(s,x),\theta}$ be a semi-infinite geodesic, let $\alpha\in(0,2/3]$ and let $f$ be a deterministic { H\"older} $\alpha^-$ function defined on a compact interval $J\subset[s,\infty)$. Then for $Z=\{t:f(t)=\gamma(t)\}$ we have $\dim
_H Z \le 1-\alpha.$
\end{lemma}
\begin{proof} By the landscape symmetries, it suffices to prove the claim for $s=x=\theta=0$. Let $I$ be a compact subinterval of $J\cap (s,\infty)$. 

By \cite[Lemma~3.11]{ganguly2022fractal}, $\gamma$ is almost surely
$\beta$-H\"older continuous on $I$ for every $\beta<2/3$. By \cite[Corollary 
7.2]{ganguly2022fractal}, there is a constant $c$ so that  for all $t\in I$,  $\gamma(t)$ has density bounded above
by $c$. The same density bound holds for $\gamma-f$. Since $f$ is $\beta$-H\"older for every $\beta<\alpha\le 2/3$, Lemma \ref{lem:zero-set-dimension} applied to $\gamma-f$ implies that $\dim_H(Z\cap I)\le 1-\alpha$. Letting $I\uparrow J$ over a countable sequence now proves the claim.
\end{proof}

\begin{lemma}[Countable covering of the geodesic network]\label{l:countable-geodesic-cover}
Almost surely, the union $\mathcal G$ of the interiors of all finite directed-landscape geodesics \eqref{e:G-def} satisfies
$$
\mathcal{G}=\big\{(\gamma_{(s,x),\theta}(r),r):x,s\in \mathbb Q, \theta \in \mathbb Z,r\in \mathbb R, s<r\big\}
$$
where $\gamma_{q,\theta}$ is the unique semi-infinite geodesic from $q$ in direction $\theta$.
\end{lemma}
\begin{proof}
This is a simple consequence of the results in \cite{rahmanVirag2025Trees}. We will use the language of that paper. 

Let $z=(\gamma(r),r)$ be an interior point of a finite geodesic $\gamma:[s,t]\to\mathbb R$. Choose $s<a<r<b<t$. The point $(\gamma(a),a)$, being interior to a geodesic, is not a forward $3$-star by Lemma 6 in \cite{rahmanVirag2025Trees}. By Proposition 12 \cite{rahmanVirag2025Trees}, the restriction $\gamma|_{[a,b]}$ is a finite concatenation of segments of semi-infinite geodesics in two directions $\theta_1<\theta_2$. Corollary 10 of that paper allows us to enlarge the direction interval, so we may take both directions to be integers. The point $z$, including when it is a switching point, is interior to at least one of those semi-infinite geodesics.

{ At} a switching point use the semi-infinite geodesic containing the preceding positive-length segment{ ; the point lies strictly after} its starting time, and belongs to { it} by continuity.
 Section 2.1 of the same paper, invoking \cite[Lemma 21]{bhatia2024duality}, says that every interior point of a geodesic in a fixed direction also lies on a geodesic in that direction started from a rational space-time point. Thus $z$ belongs to the displayed union. The reverse inclusion holds because every interior point of a semi-infinite geodesic lies in the interior of one of its finite segments.
\end{proof}

\begin{proof}[Proof of Corollary \ref{c: f_int}]  
By Lemma \ref{l:countable-geodesic-cover}
$$
\mathcal{G}=\big\{(\gamma_{(s,x),\theta}(r),r):x,s\in \mathbb Q, \theta \in \mathbb Z,r\in \mathbb R, s<r\big\}
$$
a countable union of semi-infinite geodesics. Since the Hausdorff dimension of a countable union is equal to the supremum of Hausdorff dimensions, we get  
{ \[
\dim_H\{t:(f(t),t)\in\mathcal G\}
=\sup_{\substack{x,s\in\mathbb Q\\\theta\in\mathbb Z}}
\dim_H\{t>s:f(t)=\gamma_{(s,x),\theta}(t)\}
\le\frac{57}{61}.
\]}
by Lemma \ref{l:GZ} applied with $\alpha={4/61}$. 
Finally by the same reasoning and Theorem \ref{t: main}, for any road $\pi$, 
{ \[\dim_H\{t:f(t)=\pi(t)\}\le \dim_H{ \{t:(\pi(t),t)\notin \mathcal G\}}\vee\dim_H\{t:(f(t),t)\in \mathcal G\}\le 57/61.\qedhere \]}  
\end{proof}

\subsection{Hausdorff dimension lower bound}

We now prove the lower bound stated in the introduction.

\begin{proof}[Proof of Theorem \ref{t: mainlb}]
We construct a road by joining points of a deterministic Cantor set
$\cC$ using geodesics. At step $k$ of the Cantor construction, remove
the middle proportion
\[
    \beta_k=\frac34+\frac{1}{k+10}
\]
from each surviving interval. Let $l_0=1$ and
\[
    l_k=l_{k-1}\frac{1-\beta_k}{2}.
\]
After $k$ steps, there are $2^k$ surviving intervals, each of length
$l_k$. Since
\[
    l_k
    =8^{-k}\prod_{j=1}^k\left(1-\frac{4}{j+10}\right)
    =8^{-k}k^{-4+o(1)},
\]
the usual Cantor-set dimension calculation gives
\[
    \dim_H(\cC)
    =\lim_{k\to\infty}\frac{k\log 2}{-\log l_k}
    =\frac13;
\]
see \cite[Chapter 3]{Falc14}.

Let $\cC_k$ be the endpoints of the intervals left after $k$ rounds
of deletion. Define $\pi_k$ by concatenating geodesics between
consecutive points of $\{0\}\times\cC_k$, and write
\[
    D_k=|\pi_{k-1}|-|\pi_k|.
\]
Each refinement replaces a geodesic by a path through two prescribed
intermediate points, so $D_k\geq0$.

At step $k$, we refine $2^{k-1}$ intervals of length $l_{k-1}$.
Each is divided into intervals of lengths $l_k$,
$\beta_k l_{k-1}$, and $l_k$. By stationarity and scaling of the
landscape, and since  the GUE Tracy-Widom { distribution}  of $\cL(0,0;0,1)$ has negative expectation $-a$ (\cite{tracy1994level}) we obtain
\begin{align*}
    \EE \bs{D_k}
    &=a\,2^{k-1}l_{k-1}^{1/3}
      \left[
        2\left(\frac{1-\beta_k}{2}\right)^{1/3}
        +\beta_k^{1/3}-1
      \right] \\
    &\leq C\,2^k l_{k-1}^{1/3}.
\end{align*}
Since
$
    2^k l_{k-1}^{1/3}=O(k^{-4/3}),$
we have $\sum_k\EE [D_k]<\infty$. As $D_k\geq0$, monotone convergence
implies $\sum_kD_k<\infty$ almost surely. Thus $|\pi_k|$ decreases
to a finite value $L$ a.s. 

The fixed endpoints, the lower bound $|\pi_k|\geq L$, and the
coercive landscape bound (Theorem \ref{t:DL-value-bound}) keep the
graphs of $\pi_k$ in a common compact spatial region.
Lemma \ref{l:Holder1/2} then makes the paths uniformly
equicontinuous. Pass to a uniformly convergent subsequence and call
its limit $\pi$. For any partition $0=t_0<\cdots<t_j=1$,
continuity of $\cL$ gives
\[
    \sum_{i=1}^j
      \cL\bigl(\tilde\pi(t_{i-1});\tilde\pi(t_i)\bigr)
    =
    \lim_k\sum_{i=1}^j
      \cL\bigl(\tilde\pi_k(t_{i-1});\tilde\pi_k(t_i)\bigr)
    \geq L.
\]
Taking the infimum over partitions shows that
$|\pi|\geq L>-\infty$, so $\pi$ is a road. Moreover,
$\pi(t)=0$ for every $t\in\cC$: this holds at every endpoint in
$\bigcup_k\cC_k$, and those endpoints are dense in $\cC$.

It remains to check that the non-geodesic points in $\cC$ have
dimension $1/3$. Fix $s<1/3$ and let $\mu_s$ be a Frostman
probability measure on $\cC$ with exponent $s$. Write
\[
    G=\{t\in\cC:(0,t)\in\cG\}.
\]

{ Indeed,} the countable covering lemma { and} { the} { non-atomic} { one-time} { laws} of the rational-start, integer-direction semi-infinite geodesics { show that each fixed $(0,t)$ is} almost surely { outside} { $\mathcal G$}. The countable { representation} also { makes} the { set} { in} the { following} { Fubini} { calculation} { measurable}.

\[
    \EE[\mu_s(G)]
    =\int_{\cC}\PP\bigl((0,t)\in\cG\bigr)\,d\mu_s(t)
    =0.
\]
Consequently $\mu_s(\cC\setminus G)=1$ almost surely, and the
Frostman bound implies $\dim_H(\cC\setminus G)\geq s$. Letting $s\uparrow1/3$ over  a countable sequence, we get
\[
    \dim_H\{t:(\pi(t),t)\notin\cG\}\geq\frac13 \qquad a.s. \qedhere
\]
\end{proof}

\noindent{\bf AI use.}
All proofs are due to the authors. We used ChatGPT Astra 6.1 and Claude Opus 5.5 during the final revision of the manuscript, to help find errors, correct minor mistakes, edit the exposition, and create figures. The authors take full responsibility for the contents of the paper.

\bigskip

\noindent{\bf Acknowledgements.}
B\'alint Vir\'ag was partially supported by an NSERC Discovery Grant.

\bibliography{new}

@article{SarkarSlyZhang2024,
  author  = {Sarkar, Sourav and Sly, Allan and Zhang, Lingfu},
  title   = {Infinite order phase transition in the slow bond {TASEP}},
  journal = {Communications on Pure and Applied Mathematics},
  volume  = {77},
  number  = {6},
  pages   = {3107--3140},
  year    = {2024},
  doi     = {10.1002/cpa.22185}
}

@article{dauvergne2022disjoint,
  author  = {Dauvergne, Duncan and Zhang, Lingfu},
  title   = {Disjoint optimizers and the directed landscape},
  journal = {Mem. Amer. Math. Soc.},
  volume  = {303},
  number  = {1524},
  year    = {2024},
  doi     = {10.1090/memo/1524}
}

@book{Falc14,
title = "Fractal geometry: mathematical foundations and applications",
author = "Falconer, Kenneth",
year = "2014",
month = jan,
language = "English",
isbn = "9781119942399",
publisher = "John Wiley \& Sons, Ltd",
address = "United States",
edition = "3rd",
}

@article{tracy1994level,
	title={Level-spacing distributions and the {A}iry kernel},
	author={Tracy, Craig A and Widom, Harold},
	journal={Communications in Mathematical Physics},
	volume={159},
	number={1},
	pages={151--174},
	year={1994},
	publisher={Springer}
}

@article{basu2014last,
  title={Last passage percolation with a defect line and the solution of the slow bond problem},
  author={Basu, Riddhipratim and Sidoravicius, Vladas and Sly, Allan},
  journal={arXiv:1408.3464},
  year={2014}
}

@article{dauvergne2024wiener,
  title={Wiener densities for the {Airy} line ensemble},
  author={Dauvergne, Duncan},
  journal={Proceedings of the London Mathematical Society},
  volume={129},
  number={4},
  pages={e12638},
  year={2024},
  publisher={Wiley Online Library}
}

@article{CH,
	Author = {Corwin, Ivan and Hammond, Alan},
	Journal = {Inventiones mathematicae},
	Number = {2},
	Pages = {441--508},
	Publisher = {Springer},
	Title = {{Brownian} {G}ibbs property for {A}iry line ensembles},
	Volume = {195},
	Year = {2014}}

@article{bates2019hausdorff,
  author = {Bates, Erik and Ganguly, Shirshendu and Hammond, Alan},
  title = {Hausdorff dimensions for shared endpoints of disjoint geodesics in the directed landscape},
  journal = {Electronic Journal of Probability},
  volume = {27},
  pages = {1--44},
  year = {2022},
  doi = {10.1214/21-EJP706}
}

@article{hammond2022brownian,
  author = {Hammond, Alan},
  title = {{Brownian} regularity for the {Airy} line ensemble, and multi-polymer watermelons in {Brownian} last passage percolation},
  journal = {Memoirs of the American Mathematical Society},
  volume = {277},
  number = {1363},
  year = {2022},
  doi = {10.1090/memo/1363}
}

@article{basu2019fractal,
  author = {Basu, Riddhipratim and Ganguly, Shirshendu and Hammond, Alan},
  title = {Fractal geometry of {Airy}$_2$ processes coupled via the {Airy} sheet},
  journal = {The Annals of Probability},
  volume = {49},
  number = {1},
  pages = {485--505},
  year = {2021},
  doi = {10.1214/20-AOP1444}
}

@article{busani2024stationary,
  title={The stationary horizon and semi-infinite geodesics in the directed landscape},
  author={Busani, Ofer and Sepp{\"a}l{\"a}inen, Timo and Sorensen, Evan},
  journal={The Annals of Probability},
  volume={52},
  number={1},
  pages={1--66},
  year={2024},
  publisher={Institute of Mathematical Statistics}
}

@article{rahman2021infinite,
  author = {Rahman, Mustazee and Vir{\'a}g, B{\'a}lint},
  title = {Infinite geodesics, competition interfaces and the second class particle in the scaling limit},
  journal = {Annales de l'Institut Henri Poincar{\'e}, Probabilit{\'e}s et Statistiques},
  volume = {61},
  number = {2},
  pages = {1075--1126},
  year = {2025},
  doi = {10.1214/23-AIHP1454}
}

@article{dauvergne2022three,
  title={Three-halves variation of geodesics in the directed landscape},
  author={Dauvergne, Duncan and Sarkar, Sourav and Vir{\'a}g, B{\'a}lint},
  journal={The Annals of Probability},
  volume={50},
  number={5},
  pages={1947--1985},
  year={2022},
  publisher={Institute of Mathematical Statistics}
}

@article{dauvergne202527,
  title={The 27 geodesic networks in the directed landscape},
  author={Dauvergne, Duncan},
  journal={Inventiones mathematicae},
  volume={242},
  number={1},
  pages={123--220},
  year={2025},
  publisher={Springer}
}

@article{lindvall1999,
  author = {Lindvall, Torgny},
  title = {On {Strassen}'s Theorem on Stochastic Domination},
  journal = {Electronic Communications in Probability},
  volume = {4},
  pages = {51--59},
  year = {1999},
  doi = {10.1214/ECP.v4-1005}
}

@article{strassen1965,
  author = {Strassen, V.},
  title = {The Existence of Probability Measures with Given Marginals},
  journal = {The Annals of Mathematical Statistics},
  volume = {36},
  number = {2},
  pages = {423--439},
  year = {1965},
  doi = {10.1214/aoms/1177700153}
}

@book{villani2008optimal,
  title={Optimal transport: old and new},
  author={Villani, C{\'e}dric},
  volume={338},
  year={2008},
  publisher={Springer}
}

@article{bhatia2022atypical,
  title = {Atypical stars on a directed landscape geodesic},
  author = {Bhatia, Manan},
  journal = {arXiv:2211.05734},
  year = {2022},
 }

@article{das2023law,
  title={Law of iterated logarithms and fractal properties of the {KPZ} equation},
  author={Das, Sayan and Ghosal, Promit},
  journal={The Annals of Probability},
  volume={51},
  number={3},
  pages={930--986},
  year={2023},
  publisher={Institute of Mathematical Statistics}
}

@incollection{GH24,
  title={Fractal structure in the directed landscape},
  author={Ganguly, Shirshendu and Hegde, Milind},
  booktitle={Probability and Stochastic Processes: A Volume in Honour of Rajeeva L. Karandikar},
  pages={129--147},
  year={2024},
  publisher={Springer}
}

@article{bhatia2024duality,
  title = {Duality in the Directed Landscape and Its Applications to Fractal Geometry},
  author = {Bhatia, Manan},
  journal = {International Mathematics Research Notices},
  volume = {2024},
  number = {12},
  pages = {9633--9662},
  year = {2024},
  month = {June},
  doi = {10.1093/imrn/rnae051},
}

@article{dauvergne2021bulk,
  title={BULK PROPERTIES OF THE AIRY LINE ENSEMBLE},
  author={Dauvergne, Duncan and Vir{\'a}g, B{\'a}lint},
  journal={The Annals of Probability},
  volume={49},
  number={4},
  pages={1738--1777},
  year={2021},
  publisher={JSTOR}
}

@book{morters2010brownian,
  title={{Brownian} motion},
  author={M{\"o}rters, Peter and Peres, Yuval},
  volume={30},
  year={2010},
  publisher={Cambridge University Press}
}

@article{DOV,
  author  = {Dauvergne, Duncan and Ortmann, Janosch and Vir\'ag, B\'alint},
  title   = {The directed landscape},
  journal = {Acta Math.},
  volume  = {229},
  number  = {2},
  pages   = {201--285},
  year    = {2022},
  doi     = {10.1007/s11511-022-00284-9},
  note    = {arXiv:1812.00309}
}

@article{dauvergne2021scaling,
  title={The scaling limit of the longest increasing subsequence},
  author={Dauvergne, Duncan and Vir{\'a}g, B{\'a}lint},
  journal={arXiv:2104.08210},
  year={2021}
}

@article{Pimentel2013,
  author  = {Pimentel, Leandro P. R.},
  title   = {On the location of the maximum of a continuous stochastic process},
  journal = {J. Appl. Probab.},
  volume  = {51},
  year    = {2014},
  number  = {1},
  pages   = {152--161},
  doi     = {10.1239/jap/1395771420}
}

@article{SeppalainenShen2020,
  author = {Sepp{\"a}l{\"a}inen, Timo and Shen, Xiao},
  title = {Coalescence estimates for the corner growth model with exponential weights},
  journal = {Electronic Journal of Probability},
  volume = {25},
  pages = {1--31},
  year = {2020},
  doi = {10.1214/20-EJP489}
}

@article{zhang2020optimal,
  author  = {Zhang, Lingfu},
  title   = {Optimal exponent for coalescence of finite geodesics in exponential last passage percolation},
  journal = {Electron. Commun. Probab.},
  volume  = {25},
  year    = {2020},
  pages   = {1--14},
  doi     = {10.1214/20-ECP354}
}

@inproceedings{Newman1995,
  author = {Newman, Charles M.},
  title = {A surface view of first-passage percolation},
  booktitle = {Proceedings of the International Congress of Mathematicians},
  editor = {Chatterji, S. D.},
  pages = {1017--1023},
  publisher = {Birkh{\"a}user},
  address = {Basel},
  year = {1995},
  doi = {10.1007/978-3-0348-9078-6_94}
}

@article{LG,
  author  = {Le Gall, Jean-Fran{\c c}ois},
  title   = {Geodesics in large planar maps and in the {Brownian} map},
  journal = {Acta Mathematica},
  volume  = {205},
  pages   = {287--360},
  year    = {2010},
  doi     = {10.1007/s11511-010-0056-5}
}

@article{AK,
  author  = {Angel, Omer and Kolesnik, Brett and Miermont, Gr{\'e}gory},
  title   = {Stability of geodesics in the {Brownian} map},
  journal = {The Annals of Probability},
  volume  = {45},
  number  = {5},
  pages   = {3451--3479},
  year    = {2017},
  doi     = {10.1214/16-AOP1140}
}

@article{MQ,
  author  = {Miller, Jason and Qian, Wei},
  title   = {Geodesics in the {Brownian} map: Strong confluence and geometric structure},
  journal = {Memoirs of the American Mathematical Society},
  volume  = {315},
  number  = {1602},
  year    = {2025},
  doi     = {10.1090/memo/1602}
}

@article{gwy21,
  title={Existence and uniqueness of the {Liouville} quantum gravity metric for $\gamma\in$ (0, 2)},
  author={Gwynne, Ewain and Miller, Jason},
  journal={Inventiones mathematicae},
  volume={223},
  number={1},
  pages={213--333},
  year={2021},
  publisher={Springer}
}

@article{alberts2014continuum,
  title={The continuum directed random polymer},
  author={Alberts, Tom and Khanin, Konstantin and Quastel, Jeremy},
  journal={Journal of Statistical Physics},
  volume={154},
  number={1},
  pages={305--326},
  year={2014},
  publisher={Springer}
}

@article{dem22,
  author = {Dembin, Barbara and Elboim, Dor and Peled, Ron},
  title = {Coalescence of geodesics and the {BKS} midpoint problem in planar first-passage percolation},
  journal = {Geometric and Functional Analysis},
  volume = {34},
  pages = {733--797},
  year = {2024},
  doi = {10.1007/s00039-024-00672-z}
}

@article{bus22,
  title={Universality of the geodesic tree in last passage percolation},
  author={Busani, Ofer and Ferrari, Patrik L},
  journal={The Annals of Probability},
  volume={50},
  number={1},
  pages={90--130},
  year={2022},
  publisher={Institute of Mathematical Statistics}
}

@article{ham20,
  title={Exponents governing the rarity of disjoint polymers in {Brownian} last passage percolation},
  author={Hammond, Alan},
  journal={Proceedings of the London Mathematical Society},
  volume={120},
  number={3},
  pages={370--433},
  year={2020},
  publisher={Wiley Online Library}
}

@article{ras24,
  title={Coalescence and total-variation distance of semi-infinite inverse-gamma polymers},
  author={Rassoul-Agha, Firas and Sepp{\"a}l{\"a}inen, Timo and Shen, Xiao},
  journal={Journal of the London Mathematical Society},
  volume={110},
  number={1},
  pages={e12955},
  year={2024},
  publisher={Wiley Online Library}
}

@article{ganguly2022fractal,
  title={Fractal geometry of the space-time difference profile in the directed landscape via construction of geodesic local times},
  author={Ganguly, Shirshendu and Zhang, Lingfu},
  journal={arXiv preprint arXiv:2204.01674},
  year={2022}
}

@article{BSS19,
  author  = {Basu, Riddhipratim and Sarkar, Sourav and Sly, Allan},
  title   = {Coalescence of geodesics in exactly solvable models of last passage percolation},
  journal = {Journal of Mathematical Physics},
  volume  = {60},
  number  = {9},
  pages   = {093301},
  year    = {2019},
  doi     = {10.1063/1.5093799},
  eprint  = {1704.05219},
  archiveprefix = {arXiv},
}

@article{MSZ25,
  author  = {Martin, James B. and Sly, Allan and Zhang, Lingfu},
  title   = {Convergence of the environment seen from geodesics in exponential last-passage percolation},
  journal = {Journal of the European Mathematical Society},
  volume  = {27},
  number  = {3},
  pages   = {877--970},
  year    = {2025},
  doi     = {10.4171/JEMS/1594},
  eprint  = {2106.05242},
  archiveprefix = {arXiv},
}

@book{BishopPeres2017,
  author    = {Bishop, Christopher J. and Peres, Yuval},
  title     = {Fractals in Probability and Analysis},
  series    = {Cambridge Studies in Advanced Mathematics},
  volume    = {162},
  publisher = {Cambridge University Press},
  year      = {2017},
  doi       = {10.1017/9781316460238}
}

@article{rahmanVirag2025Trees,
  author  = {Rahman, Mustazee and Vir{\'a}g, B{\'a}lint},
  title   = {The directed landscape seen from its trees},
  journal = {Electronic Communications in Probability},
  volume  = {30},
  year    = {2025},
  pages   = {1--14},
  doi     = {10.1214/25-ECP706},
  eprint  = {2410.19070},
  archivePrefix = {arXiv}
}

@article{GangulyHammondNearGround,
  author = {Ganguly, Shirshendu and Hammond, Alan},
  title = {The geometry of near ground states in {Gaussian} polymer models},
  journal = {Electronic Journal of Probability},
  volume = {28},
  pages = {1--80},
  year = {2023},
  doi = {10.1214/22-EJP895}
}

@misc{BhatiaDynamics2026,
 author={Manan Bhatia},
 title={A non-trivial dynamics on the directed landscape},
 year={2026},
 note={arXiv:2609.35124},
 url={https://arxiv.org/abs/2609.35124}
}

@misc{BhatiaTraces2026,
 author={Manan Bhatia},
 title={Geodesic traces in dynamical {Brownian} last passage percolation},
 year={2026},
 note={arXiv:2609.35103},
 url={https://arxiv.org/abs/2609.35103}
}

@misc{GangulySchweinsbergShuai2025,
 author={Shirshendu Ganguly and Jason Schweinsberg and Yubo Shuai},
 title={Branch lengths for geodesics in the directed landscape and mutation patterns in growing spatially structured populations},
 year={2025},
 note={arXiv:2506.01787},
 url={https://arxiv.org/abs/2506.01787}
}
\bibliographystyle{plainnat}

\bigskip
\begingroup
\small\raggedright
\noindent Pranay Agarwal. Department of Mathematics, University of Toronto, Toronto, Canada.
\texttt{pranay.agarwal@mail.utoronto.ca}\par
\medskip
\noindent Sourav Sarkar. Department of Pure Mathematics and Mathematical Statistics, University of Cambridge, Cambridge, United Kingdom.
\texttt{ss2871@cam.ac.uk}\par
\medskip
\noindent B\'alint Vir\'ag. Department of Mathematics, University of Toronto, Toronto, Canada.
\texttt{balint@math.toronto.edu}\par
\endgroup

\end{document}